\documentclass[11pt,a4paper]{article}
\title{\bf Stochastic differential equations with \\
coefficients in Sobolev spaces}
\author{Shizan Fang$^c$\footnote{fang@u-bourgogne.fr} ,\quad Dejun Luo$^{a,b}$,
\quad Anton Thalmaier$^a$
\vspace{2mm}\\
{\footnotesize $^a$UR Math\'{e}matiques, Universit\'{e} du
Luxembourg, 6, rue Richard Coudenhove-Kalergi,
L-1359 Luxembourg}\\
{\footnotesize $^b$Key Laboratory of Random Complex Structures and
Data Science, Academy of Mathematics and}\\
{\footnotesize Systems Science, Chinese Academy of Sciences, Beijing
100190, China} \\
{\footnotesize $^c$I.M.B, BP 47870, Universit\'e de Bourgogne,
Dijon, France} }
\date{}
\usepackage{amssymb,amsmath,amsfonts,amsthm,color,mathrsfs}

\def\R{\mathbb{R}}
\def\E{\mathbb{E}}
\def\P{\mathbb{P}}
\def\D{\mathbb{D}}
\def\Z{\mathbb{Z}}
\def\L{\mathcal{L}}
\def\d{\textup{d}}
\def\ch{{\bf 1}}

\newcommand{\ra}{\rightarrow}

\newcommand{\ee}{\varepsilon}
\def\F{\mathcal{F}}
\def\SF{\mathscr{F}}
\def\Leb{\textup{Leb}}
\def\div{\textup{div}}
\def\fin{\hfill$\square$}
\def\dis{\displaystyle}
\def\bint{-\hskip-4.1mm\int}

\def\newdot{{\kern.8pt\cdot\kern.8pt}}

\newtheorem{theorem}{Theorem}[section]
\newtheorem{lemma}[theorem]{Lemma}
\newtheorem{corollary}[theorem]{Corollary}
\newtheorem{proposition}[theorem]{Proposition}
\newtheorem{remark}[theorem]{Remark}

\newtheorem{definition}[theorem]{Definition}

\begin{document}

\maketitle
\makeatletter 
\renewcommand\theequation{\thesection.\arabic{equation}}
\@addtoreset{equation}{section}
\makeatother 

\begin{abstract}
We consider It\^o SDE $\d X_t=\sum_{j=1}^m A_j(X_t)\,\d w_t^j +
A_0(X_t)\,\d t$ on $\R^d$. The diffusion coefficients $A_1,\cdots,
A_m$ are supposed to be in the Sobolev space $W_\text{loc}^{1,p}
(\R^d)$ with $p>d$, and to have linear growth; for the drift
coefficient $A_0$, we consider two cases: (i) $A_0$ is continuous
whose distributional divergence $\delta(A_0)$ w.r.t. the Gaussian
measure $\gamma_d$ exists, (ii) $A_0$ has the Sobolev regularity
$W_\text{loc}^{1,p'}$ for some $p'>1$. Assume $\int_{\R^d}
\exp\big[\lambda_0\bigl( |\delta(A_0)| + \sum_{j=1}^m
(|\delta(A_j)|^2 +|\nabla A_j|^2)\bigr)\big] \d\gamma_d<+\infty$ for
some $\lambda_0>0$, in the case (i), if the pathwise uniqueness of
solutions holds, then the push-forward $(X_t)_\# \gamma_d$ admits a
density with respect to $\gamma_d$. In particular, if the
coefficients are bounded Lipschitz continuous, then $X_t$ leaves the
Lebesgue measure $\Leb_d$ quasi-invariant. In the case (ii), we
develop a method used by G. Crippa and C. De Lellis for ODE and
implemented by X. Zhang for SDE, to establish the existence and
uniqueness of stochastic flow of maps.
\end{abstract}



\section{Introduction}

Let $A_0, A_1, \ldots, A_m\colon\,\R^d\ra\R^d$ be continuous vector fields
on $\R^d$. We consider the following It\^{o} stochastic differential
equation on $\R^d$ (abbreviated as SDE)
  \begin{equation}\label{ItoSDE}
  \d X_t=\sum_{j=1}^m A_j(X_t)\,\d w^j_t+A_0(X_t)\,\d t, \quad X_0=x,
  \end{equation}
where $w_t=(w^1_t,\ldots, w^m_t)$ is the standard Brownian motion on
$\R^m$. It is a classical fact in the theory of SDE (see
\cite{Huang01, IkedaWatanabe89,Kunita90,RevuzYor}) that, if the
coefficients $A_j$ are globally Lipschitz continuous, then SDE
\eqref{ItoSDE} has a unique strong solution which defines a
stochastic flow of homeomorphisms on $\R^d$; however contrary to
ordinary differential equations (abbreviated as ODE), the regularity
of the homeomorphisms is only H\"older continuity of order
$0<\alpha<1$. Thus it is not clear whether the Lebesgue measure
$\Leb_d$ on $\R^d$ admits a density under the flow $X_t$. In the
case  where the vector fields $A_j, j=0,1, \ldots, m$, are in
$C_b^{\infty}(\R^d,\R^d)$, the SDE \eqref{ItoSDE} defines a flow of
diffeomorphisms, and Kunita \cite{Kunita90} showed that the measures
on $\R^d$ which have a strictly positive smooth density with respect
to $\Leb_d$ are quasi-invariant under the flow. This result was
recently generalized in \cite{Luo09} to the case where the drift
$A_0$ is allowed to be only log-Lipschitz continuous. Studies on SDE
beyond the Lipschitz setting attracted great interest during the
last years, see for instance \cite{Fang, FangImkellerZhang,
FangZhang, Krylov, KrylovRockner, LeJanRaimond1, LeJanRaimond2,
Malliavin, Zhang05a, Zhang05b}.

In the context of ODE, existence of a flow of quasi-invariant
measurable maps associated to a vector field $A_0$ belonging to
Sobolev spaces appeared first in \cite{Cruzeiro83}. In the seminar
paper \cite{DiPernaLions}, Di Perna and Lions developed transport
equations to solve ODE without involving exponential integrability
of $|\nabla A_0|$. On the other hand, L.~Ambrosio \cite{Ambrosio1}
took advantage of using continuity equations which allowed him to
construct quasi-invariant flows associated to vector fields $A_0$
with only BV regularity. In the framework for Gaussian measures, the
Di~Perna-Lions method was developed in \cite{CiprianoCruzeiro05},
also in \cite{AmbrosioFigalli09, FangLuo09} on the Wiener space.

The situation for SDE is quite different: even for the vector fields
$A_0, A_1, \ldots, A_m$ in $C^\infty$ with linear growth, if no
conditions were imposed on the growth of the derivatives, the
SDE~\eqref{ItoSDE} could not define a flow of diffeomorphisms (see
\cite{Li, LiScheutzow}). More precisely, let $\tau_x$ be the life
time of the solution to \eqref{ItoSDE} starting from $x$. The SDE
\eqref{ItoSDE} is said to be \textit{complete} if for each $x\in
\R^d$, $\P(\tau_x=+\infty)=1$; it is said to be \textit{strongly
complete} if $\P(\tau_x=+\infty,\, x\in\R^d)=1$. The goal in
\cite{LiScheutzow} is to construct examples for which the
coefficients are smooth, but the SDE \eqref{ItoSDE} is not strongly
complete (see \cite{FangImkellerZhang, Li} for positive examples).
Now consider
  $$\Sigma=\{(w,x)\in \Omega\times\R^d;\, \tau_x(w)=+\infty\}.$$
Suppose that the SDE \eqref{ItoSDE} is complete, then for any
probability measure $\mu$ on $\R^d$,
  $$\int_{\R^d}\left(\int_\Omega {\bf
  1}_{\Sigma}(w,x)\,\d\P(w)\right)\d\mu(x)=1.$$
By Fubini's theorem, $\int_{\Omega}\left(\int_{\R^d} {\bf
1}_{\Sigma}(w,x)\,\d\mu(x)\right)\d\P(w)=1$. It follows that there
exists a full measure subset $\Omega_0\subset\Omega$ such that for
all $w\in\Omega_0,\, \tau_x(w)=+\infty$ holds for $\mu$-almost every
$x\in\R^d$. Now under the existence of a complete unique strong
solution to SDE \eqref{ItoSDE}, we have a flow of measurable maps
$x\ra X_t(w,x)$.

Recently, inspired by a previous work due to Ambrosio, Lecumberry
and Maniglia \cite{ALM}, Crippa and De Lellis \cite{CrippadeLellis}
obtained some new type of estimates of perturbation for ODE whose
coefficients have Sobolev regularity. More precisely, the absence of
Lipschitz condition was filled by the following inequality: for
$f\in W_{loc}^{1,1}(\R^d)$,
  \begin{equation*}
  |f(x)-f(y)|\leq C_d\, |x-y|\, \bigl( M_R|\nabla f|(x)+M_R|\nabla
  f|(y)\bigr)
  \end{equation*}
holds for $x,y\in N^c$ and $|x-y|\leq R$, where $N$ is a negligible
set of $\R^d$ and $M_Rg$ is the maximal function defined by
  \begin{equation*}
  M_Rg(x)=\sup_{0<r\leq R}\frac{1}{\Leb_d(B(x,r))}\int_{B(x,r)}
  |g(y)|\,\d y,
  \end{equation*}
here $B(x,r)=\{ y\in\R^d;|y-x|\leq r\}$; the classical moment
estimate was replaced by estimating the quantity
  \begin{equation*}
  \int_{B(0,r)}\log\bigg(\frac{|X_t(x)-\tilde X_t(x)|}{
  \sigma}+1\bigg)\,\d x,
  \end{equation*}
where $\sigma>0$ is a small parameter. This method has recently been
successfully implemented to SDE by X. Zhang in \cite{Zhang09}.

The aim in this paper is two-fold: first we shall study absolute
continuity of the push-forward measure $(X_t)_\#\Leb_d$ with respect
to $\Leb_d$, once the SDE \eqref{ItoSDE} has a unique strong
solution; secondly we shall construct strong solutions (for almost
all initial values) using the approach mentioned above for SDE with
coefficients in Sobolev space. The key point is to obtain {\it a
priori} $L^p$ estimate for the density. To this end, we shall work
with the standard Gaussian measure $\gamma_d$; this will be done in
Section 2. The main result in Section 3 is the following

\begin{theorem}\label{th.A} Let $A_0, A_1,\ldots,A_m$ be continuous vector
fields on $\R^d$ of linear growth. Assume that the diffusion
coefficients $A_1,\ldots,A_m$ are in the Sobolev space
$\cap_{q>1}\D_1^q(\gamma_d)$ and that $\delta(A_0)$ exists;
furthermore there exists a constant $\lambda_0>0$ such that
  \begin{equation}\label{th.A-1}
  \int_{\R^d}\exp\bigg[\lambda_0\bigg(|\delta(A_0)|+\sum_{j=1}^m
  \big(|\delta(A_j)|^2+|\nabla A_j|^2\big)\bigg)\bigg]\d\gamma_d<+\infty.
  \end{equation}
Suppose that pathwise uniqueness holds for \textup{SDE}
\eqref{ItoSDE}. Then $(X_t)_\# \gamma_d$ is absolutely continuous
with respect to $\gamma_d$ and the density is in the space
$L^1\log{L^1}$.
\end{theorem}

A consequence of this theorem concerns the following classical
situation.

\begin{theorem}\label{th.B} Let $A_0, A_1, \ldots, A_m$ be globally Lipschitz continuous.
Suppose that there exists a constant $C>0$ such that
  \begin{equation}\label{th.B-1}
  \sum_{j=1}^m\langle x,A_j(x)\rangle^2\leq C\,(1+|x|^2)
  \quad \mbox{for all }x\in\R^d.
  \end{equation}
Then the stochastic flow of homeomorphisms $X_t$ generated by \textup{SDE}
\eqref{ItoSDE} leaves the Lebesgue measure $\Leb_d$ quasi-invariant.
\end{theorem}

Remark that the condition \eqref{th.B-1} not only includes the case
of bounded Lipschitz diffusion coefficients, but also, maybe more
significant, indicates the role of dispersion: the vector fields
$A_1, \cdots, A_m$ should not go radically into infinity. The
purpose of Section 4 is to find conditions that guarantee strict
positivity of the density, in the case where the existence of the
inverse flow is not known, see Theorem \ref{th4.3}.

The main result in Section 5 is

\begin{theorem}\label{th.C}
Assume that the diffusion coefficients $A_1, \cdots A_m$ belong to
the Sobolev space $\cap_{q>1}\D_1^q(\gamma_d)$ and the drift
$A_0\in\D_1^q(\gamma_d)$ for some $q>1$. Assume \eqref{th.A-1} and
that the coefficients $A_0, A_1, \cdots, A_m$ are of linear growth,
then there is a unique stochastic flow of measurable maps $X:
[0,T]\times\Omega\times\R^d\ra \R^d$, which solves \eqref{ItoSDE}
for almost all initial $x\in\R^d$ and the push-forward
$(X_t(w,\cdot))_\#\gamma_d$ admits a density with respect to
$\gamma_d$, which is in $L^1\log{L^1}$.
\end{theorem}

When the diffusion coefficients satisfy the uniform ellipticity, a
classical result due to Stroock and Varadhan \cite{StroockVaradhan}
says that if the diffusion coefficients $A_1, \cdots, A_m$ are
bounded continuous and the drift $A_0$ is bounded Borel measurable,
then the weak uniqueness holds, that is the uniqueness in law of the
diffusion. This result was strengthened by Veretennikov
\cite{Veretennikov}, saying that in fact the pathwise uniqueness
holds. When $A_0$ is not bounded, some conditions on diffusion
coefficients were needed. In the case where the diffusion matrix
$a=(a_{ij})$ is the identity, the drift $A_0$ in \eqref{ItoSDE} can
be quite singular: $A_0\in L_{loc}^p(\R^d)$ with $p>d+2$ implies
that the SDE \eqref{ItoSDE} has the pathwise uniqueness (see
Krylov-R\"ockner \cite{KrylovRockner} for a more complete study); if
the diffusion coefficients $A_1, \cdots, A_m$ are bounded
continuous, under a Sobolev condition, namely, $A_j\in W_{loc}^{1,
2(d+1)}$ for $j=1, \cdots, m$ and $A_0\in L_{loc}^{2(d+1)}(\R^d)$,
X. Zhang proved in \cite{Zhang05a} that the SDE \eqref{ItoSDE}
admits a unique strong solution. Note that even in this uniformly
non-degenerated case, if the diffusion coefficients lose the
continuity, there are counterexamples for which the weak uniqueness
does not hold, see \cite{Krylov, Safonov}.

Finally we would like to mention that under weaker Sobolev type
conditions, the connection between weak solutions and Fokker-Planck
equations was investigated in \cite{Figalli, LebrisLions}, some
notions of ``generalized solutions'', as well as the phenomena of
coalescence and splitting, were investigated in \cite{LeJanRaimond1,
LeJanRaimond2}. Stochastic transport equations were studied in
\cite{Flandoli,Zhang09}.

\section{$L^p$ estimate of the density}

The purpose of this section is to derive {\it a priori} estimates for
the density; we assume that the coefficients $A_0, A_1, \ldots,
A_m$ of SDE \eqref{ItoSDE} are {\it smooth with compact support} in
$\R^d$. Then the solution $X_t$, i.e., $x\mapsto X_t(x)$, is a stochastic flow of
diffeomorphisms on $\R^d$. Moreover SDE \eqref{ItoSDE} is equivalent
to the following Stratonovich SDE
  \begin{equation}\label{StraSDE}
  \d X_t=\sum_{j=1}^m A_j(X_t)\circ\d w^j_t+\tilde A_0(X_t)\,\d t, \quad X_0=x,
  \end{equation}
where $\tilde A_0=A_0-\frac12\sum_{j=1}^m\L_{A_j}A_j$ and $\L_A$
denotes the Lie derivative with respect to $A$.

Let $\gamma_d$ be the standard Gaussian measure on $\R^d$, and
$\gamma_t=(X_t)_\#\gamma_d,\ \tilde\gamma_t=(X_t^{-1})_\# \gamma_d$
the push-forwards of $\gamma_d$ respectively by the flow $X_t$ and
its inverse flow $X_t^{-1}$. To fix ideas, we denote by
$(\Omega,\SF,\P)$ the probability space on which the Brownian motion
$w_t$ is defined. Let $K_t=\frac{\d\gamma_t}{\d\gamma_d}$ and
$\tilde K_t= \frac{\d\tilde\gamma_t}{\d\gamma_d}$ be the densities
with respect to $\gamma_d$. By Lemma 4.3.1 in \cite{Kunita90},
the Radon-Nikodym derivative $\tilde K_t$ has the following explicit
expression
  \begin{equation}\label{Radon-Nikodym}
  \tilde K_t(x)=
  \exp\bigg(-\sum_{j=1}^m\int_0^t\delta(A_j)(X_s(x))\circ\d
  w^j_s-\int_0^t\delta(\tilde A_0)(X_{s}(x))\,\d s\bigg),
  \end{equation}
where $\delta(A_j)$ denotes the divergence of $A_j$ with respect to the
Gaussian measure $\gamma_d$:
  $$\dis \int_{\R^d} \langle \nabla\varphi, A_j\rangle\,\d\gamma_d
  =\int_{\R^d}\varphi\, \delta(A_j)\,\d\gamma_d,\quad \varphi\in
  C_c^1(\R^d).$$
It is easy to see that $K_t$ and $\tilde K_t$ are related to each
other by the equality below:
  \begin{equation}\label{relation}
  K_t(x)=\big[\tilde K_t\big(X_t^{-1}(x)\big)\big]^{-1}.
  \end{equation}
In fact, for any $\psi\in C_c^\infty(\R^d)$, we have
  \begin{eqnarray*}
  \int_{\R^d}\psi(x)\,\d\gamma_d(x)
  &=&\int_{\R^d}\psi\big[X_t\big(X_t^{-1}(x)\big)\big]\,\d\gamma_d(x)\cr
  &=&\int_{\R^d}\psi[X_t(y)]\tilde K_t(y)\,\d\gamma_d(y)=\int_{\R^d}\psi(x)\tilde
  K_t\big(X^{-1}_t(x)\big)K_t(x)\,\d\gamma_d(x),
  \end{eqnarray*}
which leads to \eqref{relation} due to the arbitrariness of $\psi\in
C_c^\infty(\R^d)$. In the following we shall estimate the
$L^p(\P\times\gamma_d)$ norm of $K_t$.

We rewrite the density \eqref{Radon-Nikodym} with the It\^{o}
integral:
  \begin{equation}\label{Radon-Nikodym-Ito}
  \tilde K_t(x)=\exp\bigg(-\sum_{j=1}^m\int_0^t\delta(A_j)(X_s(x))\,\d
  w^j_s-\int_0^t\bigg[\frac12\sum_{j=1}^m\L_{A_j}\delta(A_j) +\delta(\tilde A_0)
  \bigg](X_{s}(x))\,\d s\bigg).
  \end{equation}

\begin{lemma} We have
  \begin{equation}\label{sect-finite.3}
  \frac12\sum_{j=1}^m\L_{A_j}\delta(A_j)+\delta(\tilde A_0)
  =\delta(A_0)+\frac12\sum_{j=1}^m|A_j|^2+\frac12\sum_{j=1}^m\langle\nabla A_j,(\nabla
  A_j)^\ast\rangle,
  \end{equation}
where $\langle\cdot,\cdot\rangle$ denotes the inner product of\/
$\R^d\otimes\R^d$ and $(\nabla A_j)^\ast$ the transpose of $\nabla
A_j$.
\end{lemma}

\noindent{\bf Proof.} Let $A$ be a $C^2$ vector field on $\R^d$. From
the expression
  \begin{equation*}
  \delta(A)=\sum_{k=1}^d\bigg( x_k A^k - \frac{\partial A^k}{\partial
  x_k}\bigg),
  \end{equation*}
we get
  \begin{equation}\label{sect-finite.1}
  \L_A\delta(A)=\sum_{\ell,k=1}^d\bigg( A^\ell A^k\delta_{k\ell}+A^\ell
  x_k \frac{\partial A^k}{\partial x_\ell}-A^\ell \frac{\partial^2
  A^k}{\partial x_\ell\partial x_k}\bigg).
  \end{equation}
Note that
  \begin{equation*}
  \frac{\partial}{\partial x_k}\bigg(A^\ell \frac{\partial
  A^k}{\partial x_\ell}\bigg)= \frac{\partial A^k}{\partial
  x_\ell}\frac{\partial A^\ell}{\partial x_k} + A^\ell
  \frac{\partial^2 A^k}{\partial x_k\partial x_\ell}.
  \end{equation*}
Thus, by means of \eqref{sect-finite.1}, we obtain
  \begin{equation}\label{sect-finite.2}
  \L_A\delta(A)=|A|^2 + \delta(\L_AA) + \langle \nabla A, (\nabla
  A)^*\rangle.
  \end{equation}
Recall that $\delta(\tilde A_0)=\delta(A_0) - \frac12\sum_{j=1}^m
\delta(\L_{A_j}A_j)$. Hence, replacing $A$ by $A_j$ in
\eqref{sect-finite.2} and summing over~$j$, gives formula
\eqref{sect-finite.3}.  \fin

\medskip

We can now prove the following key estimate.

\begin{theorem}\label{sect-2-thm}
For $p>1$,
  \begin{equation}\label{densityestimate}
  \|K_t\|_{L^p(\P\times\gamma_d)}\leq\bigg[\int_{\R^d}\exp\bigg(pt\Big[2|\delta(A_0)|+\sum_{j=1}^m\big(|A_j|^2+|\nabla
  A_j|^2+2(p-1)|\delta(A_j)|^2\big)\Big]\bigg)\d\gamma_d\bigg]^{\frac{p-1}{p(2p-1)}}.
  \end{equation}
\end{theorem}

\noindent{\bf Proof.} Using relation \eqref{relation}, we have
  \begin{align}\label{sect-2-thm.1}
  \int_{\R^d}\E[K_t^p(x)]\,\d\gamma_d(x)&=\E\int_{\R^d}\big[\tilde
  K_t\big(X_t^{-1}(x)\big)\big]^{-p}\,\d\gamma_d(x)\cr
  &=\E\int_{\R^d}\big[\tilde K_t(y)\big]^{-p}\tilde K_t(y)\,\d\gamma_d(y)\cr
  &=\int_{\R^d}\E\big[\big(\tilde K_t(x)\big)^{-p+1}\big]\,\d\gamma_d(x).
  \end{align}
To simplify the notation, denote the right hand side of
\eqref{sect-finite.3} by $\Phi$. Then $\tilde K_t(x)$ rewrites as
  \begin{equation*}
  \tilde K_t(x)=\exp\bigg(-\sum_{j=1}^m\int_0^t\delta(A_j)(X_{s}(x))\,\d
  w^j_s-\int_0^t\Phi(X_{s}(x))\,\d s\bigg).
  \end{equation*}
Fixing an arbitrary $r>0$, we get
  \begin{align*}
  \big(\tilde K_t(x)\big)^{-r}&=\exp\bigg(r\sum_{j=1}^m\int_0^t\delta(A_j)(X_{s}(x))\,\d
  w^j_s+r\int_0^t\Phi(X_{s}(x))\,\d s\bigg)\cr
  &=\exp\bigg(r\sum_{j=1}^m\int_0^t\delta(A_j)(X_{s}(x))\,\d
  w^j_s-r^2\sum_{j=1}^m\int_0^t\big|\delta(A_j)(X_{s}(x))\big|^2\,\d s\bigg)\cr
  &\hskip10mm{}\times \exp\bigg(\int_0^t\Big(r^2\sum_{j=1}^m|\delta(A_j)|^2+r\Phi\Big)(X_{s}(x))\,\d
  s\bigg).
  \end{align*}
By Cauchy-Schwarz's inequality,
  \begin{eqnarray}\label{sect-finite.5}
  \E\big[\big(\tilde K_t(x)\big)^{-r}\big]&\leq&\bigg[\E\exp\bigg(2r\sum_{j=1}^m\int_0^t\delta(A_j)(X_{s}(x))\,\d
  w^j_s
  -2r^2\sum_{j=1}^m\int_0^t\big|\delta(A_j)(X_{s}(x))\big|^2\,\d s\bigg)\bigg]^{1/2}\cr
  &&{}\times\bigg[\E\exp\bigg(\int_0^t\Big(2r^2\sum_{j=1}^m|\delta(A_j)|^2+2r\Phi\Big)(X_{s}(x))\,\d
  s\bigg)\bigg]^{1/2}\cr
  &=&\bigg[\E\exp\bigg(\int_0^t\Big(2r^2\sum_{j=1}^m|\delta(A_j)|^2+2r\Phi\Big)(X_{s}(x))\,\d
  s\bigg)\bigg]^{1/2},
  \end{eqnarray}
since the first term on the right hand side of the inequality in \eqref{sect-finite.5}
is the expectation of a martingale. Let
  \begin{eqnarray*}
  \tilde\Phi_r
  &=&2r|\delta(A_0)|+r\sum_{j=1}^m\big(|A_j|^2+|\nabla
  A_j|^2+2r|\delta(A_j)|^2\big).
  \end{eqnarray*}
Then by \eqref{sect-finite.5}, along with the definition of $\Phi$ and
Cauchy-Schwarz's inequality, we obtain
  \begin{eqnarray}\label{sect-finite.6}
  \int_{\R^d}\E\big[\big(\tilde K_t(x)\big)^{-r}\big]d\gamma_d
  \leq\bigg[\int_{\R^d}\E\exp\bigg(\int_0^t\tilde\Phi_r(X_{s}(x))\,\d
  s\bigg)\,\d\gamma_d\bigg]^{1/2}.
  \end{eqnarray}
Following the idea of A.B. Cruzeiro (\cite{Cruzeiro83} Corollary
2.2, see also Theorem 7.3 in \cite{Driver97}) and by Jensen's
inequality,
  \begin{eqnarray*}
  \exp\bigg(\int_0^t\tilde\Phi_r(X_{s}(x))\,\d s\bigg)
  =\exp\bigg(\int_0^t t\,\tilde\Phi_r(X_{s}(x))\,\frac{\d s}t\bigg)
  \leq\frac1t\int_0^t e^{t\,\tilde\Phi_r(X_{s}(x))}\,\d s.
  \end{eqnarray*}
Define $I(t)=\sup_{0\leq s\leq t}\int_{\R^d}\E[K_t^p(x)]\,\d\gamma_d$.
Integrating on both sides of the above inequality and by H\"{o}lder's
inequality,
  \begin{align*}
  \int_{\R^d}\E\exp\bigg(\int_0^t\tilde\Phi_r(X_{s}(x))\,\d
  s\bigg)\d\gamma_d(x)
  &\leq\frac1t\int_0^t\E\int_{\R^d}e^{t\,\tilde\Phi_r(X_{s}(x))}\,\d\gamma_d(x)\,\d
  s\cr
  &=\frac1t\int_0^t\E\int_{\R^d}e^{t\,\tilde\Phi_r(y)}K_{s}(y)\,\d \gamma_d(y)\,\d
  s\cr&\leq\frac1t\int_0^t\big\|e^{t\,\tilde\Phi_r}\big\|_{L^q(\gamma_d)}\|K_{s}\|_{L^p(\P\times\gamma_d)}\,\d
  s\cr
  &\leq\big\|e^{t\,\tilde\Phi_r}\big\|_{L^q(\gamma_d)}\,I(t)^{1/p},
  \end{align*}
where $q$ is the conjugate number of $p$. Thus it follows from
\eqref{sect-finite.6} that
  \begin{equation}\label{sect-2-thm.2}
  \int_{\R^d}\E\big[\big(\tilde K_t(x)\big)^{-r}\big]\,\d\gamma_d(x)
  \leq\big\|e^{t\,\tilde\Phi_r}\big\|_{L^q(\gamma_d)}^{1/2}\,I(t)^{1/{2p}}.
  \end{equation}
Taking $r=p-1$ in the above estimate and by \eqref{sect-2-thm.1}, we
obtain
  $$\int_{\R^d}\E[K_t^p(x)]\,\d\gamma_d(x)\leq \big\|e^{t\,\tilde\Phi_{p-1}}\big\|_{L^q(\gamma_d)}^{1/2}\,I(t)^{1/{2p}}.$$
Thus we have $I(t)\leq\big\|e^{t\,\tilde\Phi_{p-1}}
\big\|_{L^q(\gamma_d)}^{1/2} \,I(t)^{1/{2p}}$. Solving this
inequality for $I(t)$ gives
  \begin{eqnarray*}
  \hskip-10mm&&\int_{\R^d}\E[K_t^p(x)]\,\d\gamma_d(x)\leq I(t)\leq
  \bigg[\int_{\R^d}\exp\bigg(\frac{pt}{p-1}\tilde\Phi_{p-1}(x)\bigg)\,\d\gamma_d(x)\bigg]^{\frac{p-1}{2p-1}}.
  \end{eqnarray*}
Now the desired estimate follows from the definition of
$\tilde\Phi_{p-1}$. \fin

\begin{corollary}\label{sect-2-coro}
For any $p>1$,
  \begin{equation}\label{densityestimate2}
  \|\tilde K_t\|_{L^p(\P\times\gamma_d)}\leq\left[\int_{\R^d}\exp\left((p+1)t\bigg[2|\delta(A_0)|+\sum_{j=1}^m\big(|A_j|^2+|\nabla
  A_j|^2+2p|\delta(A_j)|^2\big)\bigg]\right)\d\gamma_d\right]^{\frac{1}{2p+1}}.
  \end{equation}
\end{corollary}

\noindent{\bf Proof.} Similar to \eqref{sect-2-thm.2}, we have for
$r>0$,
  \begin{equation}\label{sect-2-coro.1}
  \int_{\R^d}\E\big[\big(\tilde K_t(x)\big)^r\big]\,\d\gamma_d(x)
  \leq\big\|e^{t\,\tilde\Phi_r}\big\|_{L^q(\gamma_d)}^{1/2}\,I(t)^{1/{2p}},
  \end{equation}
where $\tilde\Phi_r$ and $I(t)$ are defined as above. Since
$I(t)\leq \big\|e^{t\,\tilde\Phi_{p-1}}
\big\|_{L^q(\gamma_d)}^{{p}/({2p-1})}$, by taking $r=p-1$, we get
  \begin{eqnarray*}
  &&\int_{\R^d}\E\big[\big(\tilde K_t(x)\big)^{p-1}\big]\,\d\gamma_d(x)
  \leq
  \big\|e^{t\,\tilde\Phi_{p-1}}\big\|_{L^q(\gamma_d)}^{{p}/{(2p-1)}}\cr
  &&\hskip6mm=\bigg[\int_{\R^d}\exp\bigg(pt\Big[2|\delta(A_0)|+\sum_{j=1}^m\big(|A_j|^2+|\nabla
  A_j|^2+2(p-1)|\delta(A_j)|^2\big)\Big]\bigg)\d\gamma_d\bigg]^{\frac{p-1}{2p-1}}.
  \end{eqnarray*}
Replacing $p$ by $p+1$ in the last inequality gives the claimed estimate. \fin

\section{Absolute continuity under flows generated by SDEs}

Now assume that the coefficients $A_j$ in SDE \eqref{ItoSDE} are
{\it continuous} and of linear growth. Then it is well known
that SDE \eqref{ItoSDE} has a weak solution of infinite life time. In
order to apply the results of the preceding section, we shall
regularize the vector fields using the Ornstein-Uhlenbeck semigroup
$\{P_\ee\}_{\ee>0}$ on $\R^d$:
  $$P_\ee A(x)=\int_{\R^d}A\big(e^{-\ee}x+\sqrt{1-e^{-2\ee}}\,y\big)\,\d\gamma_d(y).$$
We have the following simple properties.

\begin{lemma}\label{sect-3-lem-1}
Assume that $A$ is continuous and $|A(x)|\leq C\, (1+|x|^q)$ for
some $q\geq 0$. Then
\begin{enumerate}
\item[\rm(i)] there is $C_q>0$ independent of $\ee$, such that
  $$|P_\ee A(x)|\leq C_q\,(1+|x|^q),\quad\mbox{for all }x\in\R^d;$$
\item[\rm(ii)] $P_\ee A$ converges uniformly to $A$ on any compact
subset as $\ee\ra 0$.
\end{enumerate}
\end{lemma}

\noindent{\bf Proof.} (i) Note that $\big|e^{-\ee}x+
\sqrt{1-e^{-2\ee}} \, y\big|\leq |x|+|y|$ and that there exists a
constant $C>0$ such that $(|x|+|y|)^q\leq C\, (|x|^q+|y|^q)$. Using
the growth condition on $A$, we have for some constant $C>0$
(depending on $q$),
  \begin{align*}
  |P_\ee
  A(x)|&\leq\int_{\R^d}\big|A\big(e^{-\ee}x+\sqrt{1-e^{-2\ee}}\,y\big)\big|\d\gamma_d(y)\cr
  &\leq C\,\int_{\R^d}\big(1+|x|^q +|y|^q\big)\d\gamma_d(y)
  \leq C\, \bigl(1+|x|^q+ M_q\bigr)
  \end{align*}
where $M_q=\int_{\R^d}|y|^q\,\d\gamma_d(y)$. Changing the constant
yields (i).

(ii) Fix $R>0$ and $x$ in the closed ball $B(R)$ of radius $R$, centered at $0$.
Let $R_1>R$ be arbitrary. We have
  \begin{align}\label{sect-3-lem-1.1}
  |P_\ee A(x)-A(x)|&\leq
  \int_{\R^d}\big|A\big(e^{-\ee}x+\sqrt{1-e^{-2\ee}}\,y\big)-A(x)\big|\d\gamma_d(y)\cr
  &=\bigg(\int_{B(R_1)}+\int_{B(R_1)^c}\bigg)
  \big|A\big(e^{-\ee}x+\sqrt{1-e^{-2\ee}}\,y\big)-A(x)\big|\,\d\gamma_d(y)\cr
  &=:I_1+I_2.
  \end{align}
By the growth condition on $A$, for some constant $C_q>0$,
independent of $\ee$, we have
  \begin{align*}
  I_2&\leq\int_{B(R_1)^c}\Big(\big|A\big(e^{-\ee}x+\sqrt{1-e^{-2\ee}}\,y\big)\big|+|A(x)|\Big)\,\d\gamma_d(y)\cr
  &\leq C_q\,\int_{B(R_1)^c}\big(1+R^q+|y|^q\big)\,\d\gamma_d(y),
  \end{align*}
where the last term tends to $0$ as $R_1\ra +\infty$. For given $\eta>0$,
we may take $R_1$ large enough such that $I_2<\eta$. Then there exists
$\ee_{R_1}>0$ such that for $\ee<\ee_{R_1}$ and $|y|\leq R_1$,
  $$\big|e^{-\ee}x+\sqrt{1-e^{-2\ee}}\,y\big|\leq
  e^{-\ee}R+\sqrt{1-e^{-2\ee}}\,R_1\leq R_1.$$
Note that
  $$\big|e^{-\ee}x+\sqrt{1-e^{-2\ee}}\,y-x\big|\leq \ee R + \sqrt{2\ee}
  R_1,\quad \hbox{for } |x|\leq R,\ |y|\leq R_1.$$
Since $A$ is uniformly continuous on $B(R_1)$, there exits
$\ee_0\leq \ee_{R_1}$ such that
  $$\big|A(e^{-\ee}x+\sqrt{1-e^{-2\ee}}\,y)-A(x)\big|\leq \eta\quad \hbox{for all }y\in B(R_1),\ \ee\leq \ee_0.$$
As a result, the term $I_1\leq \eta$. Therefore by
\eqref{sect-3-lem-1.1}, for any $\ee\leq \ee_0$,
  $$\sup_{|x|\leq R}|P_\ee A(x)-A(x)|\leq 2\eta.$$
The result follows from the arbitrariness of $\eta>0$. \fin

\medskip

The vector field $P_\ee A$ is smooth on $\R^d$ but does not have
compact support. We introduce cut-off functions $\varphi_\ee\in
C_c^\infty(\R^d, [0,1])$ satisfying
  $$\varphi_\ee(x)=1\ \mbox{ if }|x|\leq\frac1\ee,\quad \varphi_\ee(x)=0
  \ \mbox{ if }|x|\geq\frac1\ee+2\quad\mbox{and }\|\nabla\varphi_\ee\|_\infty \leq 1.$$
Set
  \begin{equation*}
  A^\ee_j=\varphi_\ee P_\ee A_j, \,j=0,1,\ldots, m.
  \end{equation*}
Now consider the It\^{o} SDE \eqref{ItoSDE} with $A_j$ being
replaced by $A^\ee_j\ (j=0,1,\ldots,m)$, and denote the
corresponding terms by adding the superscript $\ee$, e.g. $X^\ee_t,\,
K^\ee_t$, etc.

In the sequel, we shall give a uniform estimate to
$K^\ee_t$. To this end, we need some preparations in the spirit of
Malliavin calculus \cite{Malliavin97}. For a vector field $A$ on
$\R^d$ and $p>1$, we say that $A\in \D_1^p(\gamma_d)$ if $A\in
L^p(\gamma_d)$ and if there exists $\nabla A\colon\, \R^d\ra\R^d\otimes\R^d$
in $L^p(\gamma_d)$ such that for any $v\in\R^d$,
  \begin{equation*}
  \nabla A(x)(v)=\partial_v A:=\lim_{\eta\ra 0}\frac{A(x+\eta
  v)-A(x)}{\eta}\quad\hbox{holds in }L^{p'}(\gamma_d)\hbox{ for any }p'<p.
  \end{equation*}
For such $A\in\D_1^p(\gamma_d)$, the divergence $\delta(A)\in
L^p(\gamma_d)$ exists and the following relations hold:
  \begin{equation}\label{A1}
  \nabla P_\ee A= e^{-\ee} P_\ee(\nabla A),\quad \delta(P_\ee
  A)=e^\ee P_\ee(\delta(A)).
  \end{equation}
If $A\in L^p(\gamma_d)$, then $P_\ee A\in\D_1^p(\gamma_d)$ and
$\lim_{\ee\ra0}\|P_\ee A-A\|_{L^p}=0$.

\begin{lemma}\label{sect-3-lem-2} Assume the vector field $A\in \D_1^p(\gamma_d)$ with
$p>1$, and denote by $A^\ee=\varphi_\ee P_\ee A$. Then for
$\ee\in\,]0,1]$,
  \begin{align*}
  |\delta(A^\ee)|&\leq P_\ee\big(|A|+e|\delta(A)|\big),\cr
  |A^\ee|^2&\leq P_\ee\big(|A|^2\big),\cr
  |\nabla A^\ee|^2&\leq P_\ee\big[2\big(|A|^2+|\nabla
  A|^2\big)\big],\cr
  |\delta(A^\ee)|^2&\leq P_\ee\big[2\big(|A|^2+e^2|\delta(
  A)|^2\big)\big].
  \end{align*}
\end{lemma}

\noindent{\bf Proof.} Note that according to \eqref{A1},
$\delta(A^\ee)=\delta(\varphi_\ee P_\ee A)=\varphi_\ee e^\ee
P_\ee\delta(A)-\langle\nabla\varphi_\ee,P_\ee A\rangle$, from where the
first inequality follows. In the same way, the other results are obtained.
\fin

\medskip

Applying Theorem \ref{sect-2-thm} to $K^\ee_t$ with $p=2$, we have
  \begin{equation}\label{sect-3.1}
  \|K^\ee_t\|_{L^2(\P\times\gamma_d)}\leq\bigg[\int_{\R^d}
  \exp\bigg(2t\Big[2|\delta(A^\ee_0)|+\sum_{j=1}^m\big(|A^\ee_j|^2+|\nabla
  A^\ee_j|^2+2|\delta(A^\ee_j)|^2\big)\Big]\bigg)\d\gamma_d\bigg]^{1/6}.
  \end{equation}
By Lemma \ref{sect-3-lem-2},
  \begin{eqnarray*}
  &&2\,|\delta(A^\ee_0)|+\sum_{j=1}^m\big(|A^\ee_j|^2+|\nabla
  A^\ee_j|^2+2|\delta(A^\ee_j)|^2\big)\cr
  &&\hskip6mm\leq P_\ee\bigg[2|A_0|+2e|\delta(A_0)|+\sum_{j=1}^m\big(7|A_j|^2+2|\nabla A_j|^2+4e^2|\delta(
  A_j)|^2\big)\bigg].
  \end{eqnarray*}
We deduce from Jensen's inequality and the invariance of $\gamma_d$
under the action of the semigroup $P_\ee$ that
  \begin{equation}\label{sect-3.2}
  \|K^\ee_t\|_{L^2(\P\times\gamma_d)}\leq\bigg[\int_{\R^d}
  \exp\bigg(4t\Big[|A_0|+e|\delta(A_0)|+\sum_{j=1}^m\big(4|A_j|^2+|\nabla
  A_j|^2+2e^2|\delta(A_j)|^2\big)\Big]\bigg)\d\gamma_d\bigg]^{1/6}
  \end{equation}
for any $\ee\leq 1$. According to \eqref{sect-3.2}, we consider the
following conditions.

\vskip 2mm

\noindent{\bf Assumptions (H):}
\begin{enumerate}
\item[(A1)] For $j=1, \ldots, m$, $A_j\in \cap_{q\geq
1}\D_1^q(\gamma_d)$, $A_0$ is continuous and $\delta(A_0)$ exists.

\item[(A2)] The vector fields $A_0, A_1, \ldots, A_m$ have linear
growth.

\item[(A3)] There exists $\lambda_0>0$ such that
  \begin{equation*}
  \int_{\R^d} \exp\bigg[\lambda_0\bigg(|\delta(A_0)| + \sum_{j=1}^m
  |\delta(A_j)|^2\bigg)\bigg]\d\gamma_d<+\infty.
  \end{equation*}

\item[(A4)] There exists $\lambda_0>0$ such that
  \begin{equation*}
  \int_{\R^d} \exp\bigg(\lambda_0\sum_{j=1}^m |\nabla
  A_j|^2\bigg)\d\gamma_d<+\infty.
  \end{equation*}
\end{enumerate}

Note that by Sobolev's embedding theorem, the diffusion coefficients
$A_1, \ldots, A_m$ admit H\"{o}lder continuous versions. In what
follows, we consider these continuous versions. It is clear that
under the conditions (A2)--(A4), there exists $T_0>0$ small enough,
such that
  \begin{equation}\label{LambdaT0}
  \Lambda_{T_0}:=\bigg[\int_{\R^d}\exp\bigg(4T_0\Big[|A_0|+e|\delta(A_0)|+\sum_{j=1}^m\big(4|A_j|^2+|\nabla
  A_j|^2+2e^2|\delta(A_j)|^2\big)\Big]\bigg)\d\gamma_d\bigg]^{1/6}<\infty.
  \end{equation}
In this case, for $t\in[0,T_0]$,
  \begin{equation}\label{LambdaT0New}
  \sup_{0<\ee\leq1}\|K^\ee_t\|_{L^2(\P\times\gamma_d)}
  \leq \Lambda_{T_0}.
\end{equation}

\begin{theorem}\label{sect-3-prop-1}
Let $T>0$ be given. Under \textup{(A1)--(A4)} in Assumptions {\bf
(H)}, there are two positive constants $C_1$ and $C_2$, independent
of $\ee$, such that
  $$\sup_{0<\ee\leq 1}\E\int_{\R^d}K^\ee_t|\log K^\ee_t|\,\d\gamma_d
  \leq 2\,(C_1T)^{1/2}\Lambda_{T_0}+C_2T\Lambda_{T_0}^2,
  \quad\mbox{for all }t\in[0,T].$$
\end{theorem}

\noindent{\bf Proof.} We follow the arguments of
Proposition 4.4 in \cite{FangLuo09}. By \eqref{relation} and
\eqref{Radon-Nikodym-Ito}, we have
  \begin{equation*}
  K^\ee_t(X^\ee_t(x))=\big[\tilde K^\ee_t(x)\big]^{-1}
  =\exp\bigg(\sum_{j=1}^m\int_0^t\delta(A^\ee_j)(X^\ee_s(x))\,\d w^j_s
  +\int_0^t\Phi_\ee(X^\ee_s(x))\,\d s\bigg),
  \end{equation*}
where $$\Phi_\ee=\delta(A^\ee_0)+\frac12\sum_{j=1}^m|A^\ee_j|^2
+\frac12\sum_{j=1}^m\langle\nabla A^\ee_j,(\nabla
A^\ee_j)^\ast\rangle.$$ Thus
  \begin{align}\label{LlogL}
  &\E\int_{\R^d}K^\ee_t|\log K^\ee_t|\,\d\gamma_d=\E\int_{\R^d}\big|\log
  K^\ee_t(X^\ee_t(x))\big|\,\d\gamma_d(x)\cr
  &\hskip6mm\leq\E\int_{\R^d}\bigg|\sum_{j=1}^m\int_0^t\delta(A^\ee_j)(X^\ee_s(x))\,\d
  w^j_s\bigg|\,\d\gamma_d(x)
  +\E\int_{\R^d}\left|\int_0^t\Phi_\ee(X^\ee_s(x))\,\d s\right|\d\gamma_d(x)\cr
  &\hskip6mm=:I_1+I_2.
  \end{align}
Using Burkholder's inequality, we get
  $$\E\bigg|\sum_{j=1}^m\int_0^t\delta(A^\ee_j)(X^\ee_s(x))\,\d w^j_s\bigg|
  \leq 2\,\E\bigg[\bigg(\sum_{j=1}^m\int_0^t
  \big|\delta(A^\ee_j)(X^\ee_s(x))\big|^2\,\d s\bigg)^{1/2}\bigg].$$
For the sake of simplifying the notations, write $\Psi_\ee=\sum_{j=1}^m
|\delta(A^\ee_j)|^2$. By Cauchy's inequality,
  \begin{eqnarray}\label{LlogLI_1}
  I_1\leq2\bigg[\int_0^t\E\int_{\R^d}
  \big|\Psi_\ee(X^\ee_s(x))\big|\,\d\gamma_d(x)\d
  s\bigg]^{1/2}.
  \end{eqnarray}
Now we are going to estimate $\E\int_{\R^d} \big|\Psi_\ee(X^\ee_s(x))
\big|^{2^\alpha}\d\gamma_d(x)$ for $\alpha\in\Z_+$
which will be done inductively.
First if $s\in[0,T_0]$, then by
\eqref{sect-3.2} and \eqref{LambdaT0New}, along with Cauchy's inequality,
  \begin{align}\label{sect-3.3}
  \E\int_{\R^d} \big|\Psi_\ee(X^\ee_s(x))\big|^{2^\alpha} \,\d\gamma_d(x)
  &=\E\int_{\R^d}|\Psi_\ee(y)|^{2^\alpha}K^\ee_s(y) \,\d\gamma_d(y)\cr
  &\leq\|\Psi_\ee\|_{L^{2^{\alpha+1}}(\gamma_d)}^{2^\alpha}\|K^\ee_s\|_{L^2(\P\times\gamma_d)}\cr
  &\leq\Lambda_{T_0}\|\Psi_\ee\|_{L^{2^{\alpha+1}}(\gamma_d)}^{2^\alpha}.
  \end{align}
Now for $s\in\,]T_0,2T_0]$, we shall use the flow property of
$X^\ee_s$: let $(\theta_{T_0}w)_t:=w_{T_0+t}-w_{T_0}$ and
$X^{\ee,T_0}_t$ be the solution of the It\^{o} SDE driven by the new
Brownian motion $(\theta_{T_0}w)_t$, then
  $$X^\ee_{T_0+t}(x,w)=X^{\ee,T_0}_t\big(X^\ee_{T_0}(x,w),\theta_{T_0}w\big),
  \quad\mbox{for all }t\geq0,$$
and $X^{\ee,T_0}_t$ enjoys the same properties as $X^\ee_t$.
Therefore,
  \begin{equation*}
  \begin{split}
  \E\int_{\R^d} \big|\Psi_\ee(X^\ee_s(x))\big|^{2^\alpha} \,\d\gamma_d(x)
  &=\E\int_{\R^d} \big|\Psi_\ee\big(X^{\ee,T_0}_{s-T_0}(X^\ee_{T_0}(x))\big)\big|^{2^\alpha}
  \,\d\gamma_d(x)\\
  &=\E\int_{\R^d}\big|\Psi_\ee\big(X^{\ee,T_0}_{s-T_0}(y)\big)\big|^{2^\alpha}K^\ee_{T_0}(y)\,\d\gamma_d(y)
  \end{split}
  \end{equation*}
which is dominated, using Cauchy-Schwarz inequality
  \begin{equation*}
  \begin{split}
  &\bigg(\E\int_{\R^d}\big|\Psi_\ee\big(X^{\ee,T_0}_{s-T_0}(y)\big)\big|^{2^{\alpha+1}}\,
  \d\gamma_d(y)\bigg)^{1/2}\|K^\ee_{T_0}\|_{L^2(\P\times\gamma_d)}\\
  &\hskip6mm\leq\Big(\Lambda_{T_0}\|\Psi_\ee\|_{L^{2^{\alpha+2}}(\gamma_d)}^{2^{\alpha+1}}\Big)^{1/2}\Lambda_{T_0}
  =\Lambda_{T_0}^{1+2^{-1}}\|\Psi_\ee\|_{L^{2^{\alpha+2}}(\gamma_d)}^{2^\alpha}.
  \end{split}
  \end{equation*}
Repeating this procedure, we finally obtain, for all $s\in[0,T]$,
  \begin{equation*}
  \E\int_{\R^d} \big|\Psi_\ee(X^\ee_s(x))\big|^{2^\alpha} \,\d\gamma_d(x)
  \leq\Lambda_{T_0}^{1+2^{-1}+\ldots+2^{-N+1}}\|\Psi_\ee\|_{L^{2^{\alpha+N}}(\gamma_d)}^{2^\alpha},
  \end{equation*}
where $N\in\Z_+$ is the unique integer such that $(N-1)T_0<T\leq
NT_0$. In particular, taking $\alpha=0$ gives
  \begin{equation}\label{sect-3.4}
  \E\int_{\R^d} \big|\Psi_\ee(X^\ee_s(x))\big|\d\gamma_d(x)
  \leq\Lambda_{T_0}^2\|\Psi_\ee\|_{L^{2^{N}}(\gamma_d)}.
  \end{equation}
By Lemma \ref{sect-3-lem-2},
  $$\sup_{0<\ee\leq 1}\|\Psi_\ee\|_{L^{2^{N}}(\gamma_d)}
  \leq \bigg\|2\sum_{j=1}^m\big(|A_j|^2+e^2|\delta(
  A_j)|^2\big)\bigg\|_{L^{2^{N}}(\gamma_d)}=:C_1$$
whose right hand side is finite under the assumptions (A2)--(A4).
This along with \eqref{LlogLI_1} and \eqref{sect-3.4} leads to
  \begin{equation}\label{LlogLI_1.1}
  I_1\leq 2\,(C_1T)^{1/2}\Lambda_{T_0}.
  \end{equation}

The same manipulation works for the term $I_2$ and we get
  \begin{equation}\label{LlogLI_2}
  I_2\leq C_2T\Lambda_{T_0}^2,
  \end{equation}
where
  $$C_2=\bigg\||A_0|+e|\delta(A_0)|+\frac32\sum_{j=1}^m|A_j|^2
  +\sum_{j=1}^m|\nabla A_j|^2\bigg\|_{L^{2^{N}}(\gamma_d)}<\infty.$$
Now we draw the conclusion from \eqref{LlogL}, \eqref{LlogLI_1.1}
and \eqref{LlogLI_2}. \fin

\medskip

It follows from Theorem \ref{sect-3-prop-1} that the family
$\{K^\ee_{\textstyle\cdot}\}_{0<\ee\leq1}$ is weakly compact in
$L^1([0,T]\times\Omega\times\R^d)$. Along a subsequence,
$K^\ee_{\textstyle\cdot}$ converges weakly to some
$K_{\textstyle\cdot}^{\mathstrut}\in
L^1([0,T]\times\Omega\times\R^d)$ as $\ee\ra0$.
Let
  $$\mathcal{C}=\bigg\{u\in L^1([0,T]\times\Omega\times\R^d)\colon\
  u_t\geq0,\,\int_{\R^d}\E(u_t\log u_t)\,\d\gamma_d\leq
  2\,(C_1T)^{1/2}\Lambda_{T_0}+C_2T\Lambda_{T_0}^2\bigg\}.$$
By convexity of the function $s\ra s\log s$, it is clear that
$\mathcal{C}$ is a convex subset of $L^1([0,T]\times\Omega\times
\R^d)$. Since the weak closure of $\mathcal{C}$ coincides with the
strong one, there exists a sequence of functions $u^{(n)}\in
\mathcal{C}$ which converges to $K$ in $L^1([0,T]\times\Omega\times
\R^d)$. Along a subsequence, $u^{(n)}$ converges to $K$ almost
everywhere. Hence by Fatou's lemma, we get for almost all
$t\in[0,T]$,
  \begin{equation}\label{sect-3.5}
  \int_{\R^d}\E(K_t\log K_t)\,\d\gamma_d\leq
  2\,(C_1T)^{1/2}\Lambda_{T_0}+C_2T\Lambda_{T_0}^2.
  \end{equation}

\begin{theorem}\label{sect-3-thm-1}
Assume conditions \textup{(A1)--(A4)} and that pathwise uniqueness
holds for \textup{SDE} \eqref{ItoSDE}. Then for each $t>0$, there is
a full subset $\Omega_t\subset\Omega$ such that for all
$w\in\Omega_t$, the density $\hat K_t$ of $(X_t)_\#\gamma_d$ with
respect to $\gamma_d$ exists and $\hat K_t\in L^1\log{L^1}$.
\end{theorem}

\noindent{\bf Proof.} Under these assumptions, we can use Theorem A
in \cite{KanekoNakao}.
For the convenience of the reader, we include the statement:

\begin{theorem}[\cite{KanekoNakao}]\label{KanekoNakao}
Let $\sigma_n(x)$ and $b_n(x)$ be continuous, taking values
respectively in the space of $(d\times m)$-matrices and $\R^d$.
Suppose that
  \begin{equation*}
  \sup_n \bigl(\|\sigma_n(x)\|+|b_n(x)|\bigr)\leq C\, (1+|x|),
  \end{equation*}
and for any $R>0$,
  \begin{equation*}
  \lim_{n\ra +\infty}\sup_{|x|\leq
  R}\bigl(\|\sigma_n(x)-\sigma(x)\|+|b_n(x)-b(x)|\bigr)=0.
  \end{equation*}
Suppose further that for the same Brownian motion $B(t)$, $X_n(x,t)$
solves the \textup{SDE}
  \begin{equation*}
  \d X_n(t)=\sigma_n(X_n(t))\,\d B(t)+b_n(X_n(t))\,\d t,\quad X_n(0)=x.
  \end{equation*}
If pathwise uniqueness holds for
  \begin{equation*}
  \d X(t)=\sigma(X(t))\,\d B(t)+b(X(t))\,\d t,\quad X(0)=x,
  \end{equation*}
then for any $R>0,\, T>0$,
  \begin{equation}\label{B3}
  \lim_{n\ra +\infty}\sup_{|x|\leq R}\E\bigg(\sup_{0\leq t\leq
  T}|X_n(t,x)-X(t,x)|^2\bigg)=0.
  \end{equation}
\end{theorem}

We continue the proof of Theorem \ref{sect-3-thm-1}. By means of
Lemma \ref{sect-3-lem-1} and Theorem \ref{KanekoNakao}, for any $T,
R>0$, we get
  \begin{equation}\label{sect-3-thm-1.1}
  \lim_{\ee\ra0}\sup_{|x|\leq R}\E\bigg(\sup_{0\leq t\leq
  T}|X^\ee_t(x)-X_t(x)|^2\bigg)=0.
  \end{equation}
Now fixing arbitrary $\xi\in L^\infty(\Omega)$ and $\psi\in
C_c^\infty(\R^d)$, we have
  \begin{align}\label{sect-3-thm-1.2}
  &\E\int_{\R^d}|\xi(\newdot)|\,\big|\psi(X^\ee_t(x))-\psi(X_t(x))\big|\,\d\gamma_d(x)\cr
  &\hskip6mm\leq\|\xi\|_\infty\bigg(\int_{B(R)}+\int_{B(R)^c}\bigg)\E\left|\psi(X^\ee_t(x))-\psi(X_t(x))\right|\d\gamma_d(x)\cr
  &\hskip6mm=:J_1+J_2.
  \end{align}
By \eqref{sect-3-thm-1.1},
  \begin{align}\label{sect-3-thm-1.3}
  J_1&\leq\|\xi\|_\infty\|\nabla\psi\|_\infty\int_{B(R)}\E\left|X^\ee_t(x)-X_t(x)\right|\d\gamma_d(x)\cr
  &\leq\|\xi\|_\infty\|\nabla\psi\|_\infty\left[\sup_{|x|\leq R}\E\bigg(\sup_{0\leq t\leq
  T}|X^\ee_t(x)-X_t(x)|^2\bigg)\right]^{1/2}\ra0,
  \end{align}
as $\ee$ tends to $0$. It is obvious that
  \begin{equation}\label{sect-3-thm-1.4}
  J_2\leq 2\,\|\xi\|_\infty\,\|\psi\|_\infty\,\gamma_d(B(R)^c).
  \end{equation}
Combining \eqref{sect-3-thm-1.2}, \eqref{sect-3-thm-1.3} and
\eqref{sect-3-thm-1.4}, we obtain
  $$\limsup_{\ee\ra0}\E\int_{\R^d}|\xi|\,\big|\psi(X^\ee_t(x))-\psi(X_t(x))\big|\d\gamma_d(x)
  \leq 2\,\|\xi\|_\infty\,\|\psi\|_\infty\,\gamma_d(B(R)^c)\ra0$$
as $R\uparrow\infty$. Therefore
  \begin{equation}\label{sect-3-thm-1.5}
  \lim_{\ee\ra0}\E\int_{\R^d}\xi\,\psi(X^\ee_t(x))\,\d\gamma_d(x)
  =\E\int_{\R^d}\xi\,\psi(X_t(x))\,\d\gamma_d.
  \end{equation}

On the other hand, by Theorem \ref{sect-3-prop-1}, for each fixed
$t\in [0,T]$, up to a subsequence, $K_t^\ee$ converges weakly in
$L^1(\Omega\times\R^d)$ to some $\hat K_t$, hence
  \begin{align}\label{sect-3-thm-1.6}
  \E\int_{\R^d}\xi\,\psi(X^\ee_t(x))\,\d\gamma_d(x)
  &=\E\int_{\R^d}\xi\,\psi(y)K^\ee_t(y)\,\d\gamma_d(y)\cr
  &\ra\E\int_{\R^d}\xi\,\psi(y)\hat K_t(y)\,\d\gamma_d(y).
  \end{align}
This together with \eqref{sect-3-thm-1.5} leads to
  $$\E\int_{\R^d}\xi\,\psi(X_t(x))\,\d\gamma_d(x)
  =\E\int_{\R^d}\xi\,\psi(y)\hat K_t(y)\,\d\gamma_d(y).$$
By the arbitrariness of $\xi\in L^\infty(\Omega)$, there exists a
full measure subset $\Omega_\psi$ of $\Omega$ such that
  $$\int_{\R^d}\psi(X_t(x))\,\d\gamma_d(x)=\int_{\R^d}\psi(y)\hat K_t(y)\,\d\gamma_d(y),
  \quad\mbox{for any }\omega\in\Omega_\psi.$$
Now by the separability of $C^\infty_c(\R^d)$, there exists a full
subset $\Omega_t$ such that the above equality holds for any
$\psi\in C^\infty_c(\R^d)$. Hence $(X_t)_\#\gamma_d=\hat
K_t\gamma_d$. \fin

\begin{remark} The $K_t(w,x)$ appearing in \eqref{sect-3.5} is
defined almost everywhere. It is easy to see that $K_t(w,x)$ is a
measurable modification of $\{\hat K_t(w,x);\ t\in [0,T]\}$.
\end{remark}

\begin{remark} Beyond the Lipschitz condition, several sufficient conditions
guaranteeing pathwise uniqueness for \textup{SDE} \eqref{ItoSDE} can be
found in the literature.
For example in \textup{\cite{FangZhang}}, the authors give the
condition
  \begin{equation*}
  \sum_{j=1}^m |A_j(x)-A_j(y)|^2 \leq C\, |x-y|^2 r\big(|x-y|^2\big),\quad
  |A_0(x)-A_0(y)|\leq C\, |x-y| r\big(|x-y|^2\big),
  \end{equation*}
for $|x-y|\leq c_0$ small enough, where $r\colon\, {]0,c_0]}\ra
{]0,+\infty[}$ is $C^1$ satisfying
\begin{enumerate}
\item[\rm(i)] $\displaystyle\lim_{s\ra 0}r(s)=+\infty$,
\item[\rm(ii)] $\displaystyle\lim_{s\ra 0}\frac{sr'(s)}{r(s)}=0$, and
\item[\rm(iii)] $\displaystyle\int_0^{c_0}\frac{\d s}{sr(s)}=+\infty$.
\end{enumerate}
\end{remark}\goodbreak

\begin{corollary}\label{sect-3-coro-1}
Suppose that the vector fields $A_0, A_1, \ldots, A_m$ are globally
Lipschitz continuous and there exists a constant $C>0$, such that
  \begin{equation}\label{Lip}
  \sum_{j=1}^m\langle x,A_j(x)\rangle^2\leq C\,(1+|x|^2)
  \quad \mbox{for all }x\in\R^d.
  \end{equation}
Then $(X_t)_\#\Leb_d\ll\Leb_d$ for any $t\in[0,T]$ .
\end{corollary}

\noindent{\bf Proof.} It is obvious that hypotheses (A1),~(A2)
and (A4) are satisfied, and that for some constant $C>0$,
  $$|\delta(A_0)|(x)\leq C(1+|x|^2).$$
Hence there exists $\lambda_0>0$ such that $\int_{\R^d}\exp\left(\lambda_0
|\delta(A_0)|\right)\,\d\gamma_d<+\infty$. Finally we have
  $$\sum_{j=1}^m|\,\delta(A_j)|^2(x)
  \leq 2\sum_{j=1}^m\langle x,A_j(x)\rangle^2+2\sum_{j=1}^m\hbox{Lip}(A_j)^2.$$
Therefore, under condition \eqref{Lip}, there exists $\lambda_0>0$ such that
  $$\int_{\R^d}\exp\bigg(\lambda_0\sum_{j=1}^m|\delta(A_j)|^2\bigg)\d\gamma_d<+\infty.$$
Hence, hypothesis (A3) is satisfied as well. By Theorem
\ref{sect-3-thm-1}, we have $(X_t)_\#\gamma_d=\hat K_t\gamma_d$. Let $A$ be
a Borel subset of $\R^d$ such that $\Leb_d(A)=0$, then
$\gamma_d(A)=0$; therefore $\int_{\R^d}{\bf 1}_{\{X_t(x)\in
A\}}\,\d\gamma_d(x)=0$. It follows that
${\bf 1}_{\{X_t(x)\in A\}}=0$ for $\Leb_d$ almost every $x$,
which implies $\Leb_d(X_t\in A)=0$; this means that $(X_t)_\#\Leb_d$ is absolutely continuous
with respect to $\Leb_d$. \fin

\medskip

In the next section, we shall prove that under the conditions of
Corollary \ref{sect-3-coro-1}, the density of $(X_t)_\#\Leb_d$ with
respect to $\Leb_d$ is strictly positive, in other words, $\Leb_d$
is quasi-invariant under $X_t$.

\begin{corollary}\label{coro2}
Assume that conditions \textup{(A1)--(A4)} hold. Let $\sigma=(A_j^i)$ and suppose
that for some $C>0$,
  \begin{equation*}
  \sigma(x)\sigma(x)^*\geq C\,\textup{Id},\quad \mbox{for all
  }x\in\R^d.
  \end{equation*}
Then $(X_t)_\#\gamma_d$ is absolutely continuous with respect to
$\gamma_d$.
\end{corollary}

\noindent{\bf Proof.} The conditions (A1)--(A4) are stronger than
those in Theorem 1.1 of \cite{Zhang05a} given by X.~Zhang, so the
pathwise uniqueness holds. Hence Theorem \ref{sect-3-thm-1} applies
to this case. \fin

\section{Quasi-invariance under stochastic flow}

In the sequel, by quasi-invariance we mean that the Radon-Nikodym derivative
of the corresponding push-forward measure is strictly positive.
First we prove that in the situation of
Corollary \ref{sect-3-coro-1}, the Lebesgue measure is in fact
quasi-invariant under the stochastic flow of homeomorphisms. To this
end, we need some preparations. In what follows, $T_0>0$ is chosen small
enough such that \eqref{LambdaT0} holds.

\begin{proposition}\label{prop4.1}
Let $q\geq 2$. Then
  \begin{equation}\label{4.1}
  \lim_{\ee\ra0} \int_{\R^d}\E\bigg(\bigg|\sup_{0\leq t\leq T_0}\sum_{j=1}^m \int_0^t
  \bigl[\delta(A_j^\ee)(X_s^\ee)-\delta(A_j)(X_s)\bigr]\d w_s^j\bigg|^q\bigg)\d\gamma_d=0.
  \end{equation}
\end{proposition}

\noindent{\bf Proof.} By Burkholder's inequality,
  \begin{align*}
  &\E\bigg(\sup_{0\leq t\leq T_0}\bigg|\sum_{j=1}^m \int_0^t
  \big[\delta(A_j^\ee)(X_s^\ee)-\delta(A_j)(X_s)\big]\,\d
  w_s^j\bigg|^q\bigg)\cr
  &\hskip6mm\leq C\, \E\bigg[\bigg(\int_0^{T_0} \sum_{j=1}^m
  |\delta(A_j^\ee)(X_s^\ee)-\delta(A_j)(X_s)|^2\,\d
  s\bigg)^{q/2}\bigg]\cr
  &\hskip6mm\leq C\, T_0^{q/2-1}\sum_{j=1}^m\int_0^{T_0}
  \E\bigl(|\delta(A_j^\ee)(X_s^\ee)-\delta(A_j)(X_s)|^q\bigr)\,\d s.
  \end{align*}
Again by the inequality $(a+b)^q\leq C_q\,(a^q + b^q)$, there exists
a constant $C_{q,T_0}>0$ such that the above quantity is dominated
by
  \begin{equation}\label{4.2}
  C_{q,T_0}\,\sum_{j=1}^m
  \bigg[\int_0^{T_0}\E\big(|\delta(A_j^\ee)(X_s^\ee)-\delta(A_j)(X_s^\ee)|^q\big)\d s
  +\int_0^{T_0}\E\big(|\delta(A_j)(X_s^\ee)-\delta(A_j)(X_s)|^q\big)\d s\bigg].
  \end{equation}
Let $I_1^\ee$ and $I_2^\ee$ be the two terms in the squared bracket
of \eqref{4.2}. Note that
  \begin{align}\label{4.3}
  \int_{\R^d}\E&\big(|\delta(A_j^\ee)(X_s^\ee)-\delta(A_j)(X_s^\ee)|^q\big)\d\gamma_d\cr
  &=\E\int_{\R^d} |\delta(A_j^\ee)-\delta(A_j)|^q K_s^\ee\,\d\gamma_d\cr
  &\leq\|\delta(A_j^\ee)-\delta(A_j)\|_{L^{2q}(\gamma_d)}^{q}\|K_s^\ee\|_{L^2(\P\times\gamma_d)}.
  \end{align}
According to \eqref{LambdaT0}, for $s\leq T_0$, we have
$\|K_s^\ee\|_{L^2(\P\times\gamma_d)}\leq \Lambda_{T_0}$. Remark that
  $$\delta(A_j^\ee)=\delta(\varphi_\ee P_\ee A_j)
  =\varphi_\ee\,e^\ee P_\ee\delta(A_j)-\langle\nabla\varphi_\ee,P_\ee
  A_j\rangle,$$
which converges to $\delta(A_j)$ in $L^{2q}(\gamma_d)$. By
\eqref{4.3},
  \begin{align*}
  \int_{\R^d} I_1^\ee\,\d\gamma_d
  &=\int_0^{T_0}\bigg[\int_{\R^d}\E\big(|\delta(A_j^\ee)(X_s^\ee)-\delta(A_j)(X_s^\ee)|^q\big)\d\gamma_d\bigg]\d
  s\cr
  &\leq T_0 \Lambda_{T_0}\,\|\delta(A_j^\ee)-\delta(A_j)\|_{L^{2q}(\gamma_d)}^{q}
  \end{align*}
which tends to $0$ as $\ee\ra 0$.

For the estimate of $I_2^\ee$, we remark that
$\int_{\R^d}|\delta(A_j)|^{2q}\,\d\gamma_d<+\infty$. Let $\eta>0$ be
given. There exists $\psi\in C_c(\R^d)$ such that
  $$\dis \int_{\R^d}|\delta(A_j)-\psi|^{2q}\,\d\gamma_d\leq \eta^2.$$
We have, for some constant $C_q>0$,
  \begin{align}\label{4.4}
  &\int_{\R^d}\E\big(|\delta(A_j)(X_s^\ee)-\delta(A_j)(X_s)|^q\big)\,\d\gamma_d\cr
  &\hskip6mm\leq C_q\bigg[\int_{\R^d}\E\big(|\delta(A_j)(X_s^\ee)-\psi(X_s^\ee)|^q\big)\d\gamma_d
  +\int_{\R^d}\E\big(|\psi(X_s^\ee)-\psi(X_s)|^q\big)\d\gamma_d\cr
  &\hskip18mm+\int_{\R^d}\E\big(|\psi(X_s)-\delta(A_j)(X_s)|^q\big)\d\gamma_d\bigg].
  \end{align}
Again by \eqref{LambdaT0New}, we find
  \begin{align*}
  \E\left[\int_{\R^d}|\delta(A_j)(X_s^\ee)-\psi(X_s^\ee)|^q\,\d\gamma_d\right]
  &=\E\left[\int_{\R^d}|\delta(A_j)-\psi|^q\,K_s^\ee\,\d\gamma_d\right]\cr
  &\leq\|\delta(A_j)-\psi\|_{L^{2q}(\gamma_d)}^q\,\Lambda_{T_0}\leq
  \Lambda_{T_0}\eta.
  \end{align*}
In the same way,
  \begin{equation*}
  \E\bigg[\int_{\R^d}|\delta(A_j)(X_s)-\psi(X_s)|^q\,\d\gamma_d\bigg]\leq
  \Lambda_{T_0}\eta.
  \end{equation*}
To estimate the second term on the right hand side of \eqref{4.4},
we use Theorem \ref{KanekoNakao}: from \eqref{B3}, we see that up to
a subsequence, $X_s^\ee(w,x)$ converges to $X_s(w,x)$, for each
$s\leq T_0$ and almost all $(w,x)\in\Omega\times\R^d$. By Lebesgue's
dominated convergence theorem,
  \begin{equation*}
  \lim_{\ee\ra0}\int_{\R^d}\E\big(|\psi(X_s^\ee)-\psi(X_s)|^q\big)\,\d\gamma_d=0.
  \end{equation*}
In conclusion, $\lim\limits_{\ee\ra 0}\int_{\R^d} I_2^\ee\,\d\gamma_d=0$.
According to \eqref{4.2}, the proof of \eqref{4.1} is complete. \fin

\begin{proposition}\label{prop4.2} Let $\Phi$ be defined by
  \begin{equation}\label{4.5}
  \Phi=\delta(A_0)+\frac12\sum_{j=1}^m|A_j|^2+\frac12\sum_{j=1}^m\langle\nabla
  A_j,(\nabla A_j)^\ast\rangle,
  \end{equation}
and analogously $\Phi_\ee$ where $A_j$ is replaced by $A_j^\ee$. Then
  \begin{equation}\label{4.6}
  \lim_{\ee\ra
  0}\int_{\R^d}\!\int_0^{T_0}\E\big(|\Phi_\ee(X_s^\ee)-\Phi(X_s)|^q\big)\,\d s\d\gamma_d=0.
  \end{equation}
\end{proposition}

\noindent{\bf Proof.} Along the lines of the proof of Proposition
\ref{prop4.1}, it is sufficient to remark that
  \begin{equation}\label{4.7}
  \lim_{\ee\ra 0} \|\Phi_\ee-\Phi\|_{L^{2q}(\gamma_d)}=0.
  \end{equation}
To see this, let us check convergence for the last term in the
definition of $\Phi_\ee$. We have
  \begin{eqnarray*}
  &&|\langle\nabla A_j^\ee,(\nabla A_j^\ee)^\ast\rangle-\langle\nabla
  A_j,(\nabla A_j)^\ast\rangle|\\
  &&\hskip6mm\leq \|\nabla A_j^\ee-\nabla A_j\|\,\|\nabla A_j^\ee\| + \|\nabla
  A_j\|\, \|\nabla A_j^\ee-\nabla A_j\|.
  \end{eqnarray*}
Note that $A_j^\ee=\varphi_\ee P_\ee A_j$. Thus
  $$\dis \nabla A_j^\ee = \nabla\varphi_\ee\otimes P_\ee A_j
  + e^{-\ee}\varphi_\ee P_\ee(\nabla A_j),$$
which converges to $\nabla A_j$ in $L^{2q}(\gamma_d)$ as $\ee\ra 0$.
\fin

\medskip

Now we can prove

\begin{proposition}
Under the conditions of Corollary \textup{\ref{sect-3-coro-1}}, the Lebesgue
measure $\Leb_d$ is quasi-invariant under the stochastic flow.
\end{proposition}

\noindent{\bf Proof.} Let $k_t$ be the density of $(X_t)_\#\Leb_d$
with respect to $\Leb_d$. We shall prove that $k_t$ is strictly
positive. Set
  \begin{equation}\label{A2}
  \tilde K_t^\ee(x) =\exp\bigg(-\sum_{j=1}^m
  \int_0^t\delta(A_j^\ee)(X_s^\ee(x))\,\d w_s^j -\int_0^t
  \Phi_\ee(X_s^\ee(x))\,\d s\bigg),
  \end{equation}
where $\Phi_\ee$ is defined in Proposition \ref{prop4.2}. By
\eqref{relation} we have
  \begin{equation}\label{A4}
  \int_{\R^d}\psi(X_t^\ee)\tilde
  K_t^\ee\,\d\gamma_d=\int_{\R^d}\psi\,\d\gamma_d, \quad \psi\in
  C_c^1(\R^d).
  \end{equation}
Applying Propositions \ref{prop4.1} and \ref{prop4.2}, up to a
subsequence, for each $t\leq T_0$ and almost every $(w,x)$, the term
$\tilde K_t^\ee(w,x)$ defined in \eqref{A2} converges to
  \begin{equation}\label{4.8}
  \tilde K_t(x)=\exp\bigg(-\sum_{j=1}^m \int_0^t
  \delta(A_j)(X_s(x))\,\d w_s^j - \int_0^t \Phi(X_s(x))\,\d s\bigg).
  \end{equation}
By Corollary \ref{sect-2-coro} and Lemma \ref{sect-3-lem-2}, we may
assume that $T_0$ is small enough so that for any $t\leq T_0$, the
family $\{\tilde K_t^\ee:\ee\leq 1\}$ is also bounded in
$L^2(\P\times \gamma_d)$. Therefore, by the uniform integrability,
letting $\ee\ra 0$ in \eqref{A4}, we get $\P$-almost surely,
  \begin{equation}\label{sect-4-prop-3}
  \int_{\R^d}\psi(X_t)\tilde
  K_t\,\d\gamma_d=\int_{\R^d}\psi\,\d\gamma_d, \quad \psi\in
  C_c^1(\R^d).
  \end{equation}

Now taking a Borel version of $x\ra\tilde K_t(w,x)$. Under the
assumptions, the solution $X_t$ is a stochastic flow of
homeomorphisms, hence the inverse flow $X^{-1}_t$ exists.
Consequently, if $t\leq T_0$, we deduce from \eqref{sect-4-prop-3}
that the density $K_t(w,x)$ of $(X_t)_\#\gamma_d$ with respect to
$\gamma_d$ admits the expression $K_t(w,x)=\big[\tilde
K_t\big(w,X_t^{-1}(w,x)\big)\big]^{-1}$ which is strictly positive.
For $X_{t+T_0}$ with $t\leq T_0$, we use the flow property:
$X_{t+T_0}(w,x)=X_t(\theta_{T_0}w,X_{T_0}(w,x))$. Thus, for any
$\psi\in C_c^1(\R^d)$,
  \begin{align*}
  \int_{\R^d}\psi(X_{t+T_0})\,\d\gamma_d&=\int_{\R^d}\psi\big(X_t(X_{T_0})\big)\d\gamma_d\cr
  &=\int_{\R^d}\psi(X_t)K_{T_0}\,\d\gamma_d=\int_{\R^d}\psi\, K_{T_0}\big(X_t^{-1}\big)K_t\,\d\gamma_d.
  \end{align*}
That is to say, the density $K_{t+T_0}=K_{T_0}\big(X_t^{-1}\big)
K_t$ is strictly positive. Continuing in this way, we obtain that $K_t$ is strictly positive for
any $t\geq 0$.

Now if $\rho(x)$ denotes the density of $\gamma_d$ with respect to
$\Leb_d$, then
  $$ k_t(w,x)=\rho\big(X_t^{-1}(w,x)\big)^{-1}K_t(w,x)\rho(x)>0$$
which concludes the proof. \fin

\medskip

In what follows, we will give examples for which existence of
the inverse flow is not known.

\begin{theorem}\label{th4.3} Let $A_1,\ldots,A_m$ be bounded $C^1$ vector
fields on $\R^d$ such that their derivatives are of linear growth; furthermore let
$A_0$ be continuous of linear growth such that $\delta(A_0)$
exists. Define
  \begin{equation}\label{4.9}
  \hat A_0=A_0 - \sum_{j=1}^m \L_{A_j}A_j.
  \end{equation}
Suppose that $\delta(\hat A_0)$ exists and that
  \begin{equation}\label{4.10}
  \int_{\R^d} \exp\big(\lambda_0 \,(|\delta(A_0)|+|\delta(\hat
  A_0)|)\big)\,\d\gamma_d<+\infty,\quad \hbox{for some }\lambda_0>0.
  \end{equation}
If pathwise uniqueness holds both for \textup{SDE} \eqref{ItoSDE} and for
  \begin{equation}\label{4.11}
  \d Y_t=\sum_{j=1}^m A_j(Y_t)\,\d w_t^j - \hat A_0(Y_t) \,\d t,
  \end{equation}
then the solution $X_t$ to \textup{SDE} \eqref{ItoSDE} leaves the Gaussian
measure $\gamma_d$ quasi-invariant.
\end{theorem}

\noindent{\bf Proof.} Obviously the conditions in Theorem
\ref{sect-3-thm-1} are satisfied; hence $(X_t)_\#\gamma_d
=K_t\,\gamma_d$. Let $t>0$ be given, we consider the dual SDE to
\eqref{ItoSDE}:
  \begin{equation*}
  \d Y_s^t=\sum_{j=1}^m A_j(Y_s^t)\,\d w_s^{t,j} - \hat A_0(Y_s^t)\,\d s
  \end{equation*}
for which pathwise uniqueness holds; here $w_s^t = w_{t-s}-w_t$ with
$s\in[0,t]$. Let $A_j^\ee$, $j=0,1, \ldots, m$, be the vector fields
defined as above. Consider
  \begin{equation*}
  \d Y_s^{t,\ee}=\sum_{j=1}^m A_j^\ee(Y_s^{t,\ee})\,\d w_s^{t,j} - \hat
  A_0^\ee(Y_s^{t,\ee}) \,\d s,
  \end{equation*}
where $\hat A_0^\ee = A_0^\ee - \sum_{j=1}^m \L_{A_j^\ee}A_j^\ee$.
Then it is known that $(X_t^\ee)^{-1}=Y_t^{t,\ee}$. It is easy to
check that for some constant $C>0$ independent of $\ee$,
  \begin{equation}\label{4.12}
  |\hat A_0^\ee(x)|\leq C\, (1+|x|).
  \end{equation}
Moreover,
  \begin{equation*}
  \L_{A_j^\ee}A_j^\ee =\sum_{k=1}^d
  (A_j^\ee)^k\bigg[\frac{\partial\varphi_\ee}{\partial x_k}\, P_\ee
  A_j +\varphi_\ee e^{-\ee}P_\ee\bigg(\frac{\partial A_j}{\partial
  x_k}\bigg)\bigg]
  \end{equation*}
which converges locally uniformly to $\L_{A_j}A_j$. Therefore $\hat
A_0^\ee$ converges uniformly over any compact subset to $\hat A_0$.
By Theorem \ref{KanekoNakao},
  \begin{equation*}
  \lim_{\ee\ra 0}\sup_{|x|\leq R}\E\bigg(\sup_{0\leq s\leq
  t}\big|Y_s^{t,\ee}-Y_s^t\big|^2\bigg)=0.
  \end{equation*}
It follows that, along a sequence, $Y_t^{t,\ee}$ converges to
$Y_t^t$ for almost every $(w,x)$. Now let $\psi_1, \psi_2\in
C_b(\R^d)$, we have for $t\leq T_0$,
  \begin{equation*}
  \int_{\R^d}\psi_1\cdot\psi_2(X_t^\ee)\tilde K_t^\ee\,\d\gamma_d
  =\int_{\R^d}\psi_1\big(Y_t^{t,\ee}\big)\cdot\psi_2\,\d\gamma_d.
  \end{equation*}
Letting $\ee\ra 0$ leads to
  \begin{equation}\label{4.13}
  \int_{\R^d}\psi_1\cdot\psi_2(X_t)\tilde K_t\,\d\gamma_d
  =\int_{\R^d}\psi_1\big(Y_t^{t}\big)\cdot\psi_2\, \d\gamma_d.
  \end{equation}

Taking $\psi_1$ and $\psi_2$ positive in \eqref{4.13} and using
a monotone class argument, we see that equation \eqref{4.13} holds for any
positive Borel functions $\psi_1$ and $\psi_2$. Hence taking a Borel
version of $\tilde K_t$ and setting $\psi_1=1/\tilde K_t$ in
\eqref{4.13}, we get
  \begin{equation}\label{4.14}
  \int_{\R^d}\psi_2(X_t)\,\d\gamma_d=\int_{\R^d} \big[\tilde
  K_t(Y_t^t)\big]^{-1}\psi_2\,\d\gamma_d.
  \end{equation}
It follows that $K_t=\big[\tilde K_t(Y_t^t)\big]^{-1}>0$ for $t\leq
T_0$. For $X_{t+T_0}$ with $t\leq T_0$, we shall use repeatedly
\eqref{4.13}. By the flow property, $X_{t+T_0}(w,x)=
X_t(\theta_{T_0}w, X_{T_0}(w,x))$ where $(\theta_{T_0}w)_t=
w_{t+T_0}-w_{T_0}$. Letting $t=T_0$ and replacing $\psi_2$ by
$\psi_2(X_t)$ we get
  \begin{equation*}
  \int_{\R^d} \psi_1\cdot\psi_2(X_{t+T_0})\tilde K_{T_0}\,\d\gamma_d
  =\int_{\R^d}\psi_1\big(Y_{T_0}^{T_0}\big)\psi_2(X_t)\,\d\gamma_d.
  \end{equation*}
Taking $\psi_1=1/\tilde K_{T_0}$ in the above equality, we get
  \begin{align*}
  \int_{\R^d} \psi_2(X_{t+T_0})\,\d\gamma_d
  &=\int_{\R^d}\big[\tilde
  K_{T_0}\big(Y_{T_0}^{T_0}\big)\big]^{-1}\,\psi_2(X_t)\,\d\gamma_d\cr
  &=\int_{\R^d} \big[\tilde
  K_{T_0}\big(Y_{T_0}^{T_0}\big)\big]^{-1}\,\psi_2(X_t)\tilde
  K_t^{-1}\tilde K_t\,\d\gamma_d\cr
  &=\int_{\R^d} \big[\tilde
  K_{T_0}\big(Y_{T_0}^{T_0}(Y_t^t)\big)\big]^{-1}\big[\tilde
  K_t(Y_t^t)\big]^{-1}\,\psi_2\,\d\gamma_d,
  \end{align*}
where in the last equality we have used \eqref{4.13} with
$\psi_1=\big[\tilde K_{T_0}\big(Y_{T_0}^{T_0}\big)\big]^{-1}\tilde
K_t^{-1}$. It follows that the density $K_{t+T_0}$ of
$(X_{t+T_0})_\# \gamma_d$ with respect to $\gamma_d$ is strictly
positive, and so on. \fin

\begin{corollary}
Let $A_1, \ldots, A_m$ be bounded $C^2$ vector fields such that
their derivatives up to order~$2$ grow at most linearly, and let $A_0$ be a
continuous vector field of linear growth. Suppose that
  \begin{equation}\label{Osgood}
  |A_0(x)-A_0(y)|\leq C_R\, |x-y|\log_k\frac{1}{|x-y|}\quad\hbox{for }
  |x|\leq R,\ |y|\leq R,\ |x-y|\leq c_0\ \hbox{small enough,}
  \end{equation}
where $\log_k
s=(\log{s})(\log{\log{s}})\ldots(\log{\ldots{\log{s}}})$. Suppose further that
$$\div(A_0)=\sum_{j=1}^d \frac{\partial A_0^j}{\partial x_j}$$ exists
and is bounded. Then the stochastic flow $X_t$ defined by
\textup{SDE} \eqref{ItoSDE} leaves the Lebesgue measure
quasi-invariant.
\end{corollary}

\noindent{\bf Proof.} It is obvious that $\hat A_0$ defined in
\eqref{4.9} satisfies condition \eqref{Osgood}; therefore by
\cite{FangZhang}, pathwise uniqueness holds for SDE \eqref{ItoSDE}
and \eqref{4.11}. Note that $\delta(A_0)=\langle x,A_0\rangle -
\div(A_0)$. Then condition \eqref{4.10} is satisfied; thus Theorem
\ref{th4.3} yields the result. \fin

\section{The case $A_0$ in Sobolev spaces}

From now on, $A_0$ is not supposed to be continuous, but in some
Sobolev space, that is, we replace the condition (A1) in {\bf (H)}
by
  \begin{enumerate}
  \item[(A1')] For $i=1, \ldots, m$, $A_i\in \cap_{q\geq
  1}\D_1^q(\gamma_d)$, $A_0\in\D_1^q(\gamma_d)$ for some $q>1$.
  \end{enumerate}
First we establish the following {\it a priori} estimate on
perturbations, using the method developed in \cite{Zhang09}. Let
$\{A_0, A_1, \cdots, A_m\}$ be a family of measurable vector fields
on $\R^d$. We shall give a precise definition of solution to the
following SDE
  \begin{equation}\label{5.1}
  \d X_t =\sum_{i=1}^m A_i(X_t)\,\d w_t^i + A_0(X_t)\,\d t,\quad
  X_0=x.
  \end{equation}

\begin{definition}\label{def.5.1}
We say that a measurable map $X: \Omega\times\R^d\ra C([0,T],\R^d)$
is a solution to It\^{o} SDE \eqref{5.1} if
\begin{enumerate}
\item[\rm(i)] for each $t\in [0,T]$ and almost all $x\in\R^d$, $w\ra
X_t(w,x)$ is measurable with respect to $\F_t$, i.e. the natural
filtration generated by the Brownian motion $\{w_s;\ s\leq t\}$;

\item[\rm(ii)] for each $t\in [0,T]$, there exists $K_t\in
L^1(\P\times\R^d)$ such that $(X_t(w,\cdot))_\#\gamma_d$ admits
$K_t$ as the density with respect to $\gamma_d$;

\item[\rm(iii)] almost surely
  \begin{equation*}
  \sum_{i=1}^m\int_0^T|A_i(X_s(w,x))|^2\,\d s+\int_0^T|A_0(X_s(w,x))|\,\d
  s<+\infty;
  \end{equation*}
\item[\rm(iv)] for almost all $x\in\R^d$,
  \begin{equation*}
  X_t(w,x)=x+\sum_{i=1}^m \int_0^t A_i(X_s(w,x))\,\d w_s^i+\int_0^t
  A_0(X_s(w,x))\,\d s;
  \end{equation*}
\item[\rm(v)] the flow property holds
  \begin{equation*}
  X_{t+s}(w,x)=X_t(\theta_sw, X_s(w,x)).
  \end{equation*}
\end{enumerate}
\end{definition}

Now consider another family of measurable vector fields $\{\hat A_0,
\hat A_1, \cdots, \hat A_m\}$ on $\R^d$, and denote by $\hat X_t$
the solution to the SDE
  \begin{equation}\label{5.2}
  \d\hat X_t =\sum_{i=1}^m \hat A_i(\hat X_t)\,\d w_t^i + \hat A_0(\hat
  X_t)\,\d t,\quad \hat X_0=x.
  \end{equation}
Let $\hat K_t$  be the density of $(\hat X_t)_\#\gamma_d$ and define
  \begin{equation}\label{2.3}
  \Lambda_{p,T}=\sup_{0\leq t\leq T}\Bigl(
  \|K_t\|_{L^p(\P\times\gamma_d)}\vee \|\hat
  K_t\|_{L^p(\P\times\gamma_d)}\Bigr).
  \end{equation}

\begin{theorem}\label{th.5.1}
Let $q>1$. Suppose that $A_1, \cdots, A_m$ as well as $\hat A_1,
\cdots, \hat A_m$ are in $\D_1^{2q}(\gamma_d)$ and $A_0, \hat
A_0\in\D_1^q(\gamma_d)$. Then for any $T>0$ and $R>0$, there exist
constants $C_{d,q,R}>0$ and $C_T>0$ such that for any $\sigma>0$,
  \begin{equation*}
  \begin{split}
  &\E\biggl[\int_{G_R}\log\bigg(\frac{\sup_{0\leq t\leq T}|X_t-\hat
  X_t|^2}{\sigma^2}+1\bigg)\d\gamma_d\biggr]\\
  &\hskip6mm\leq C_T\Lambda_{p,T}\bigg\{
  C_{d,q,R}\bigg[\|\nabla A_0\|_{L^q}+\bigg(\sum_{i=1}^m\|\nabla A_i\|_{L^{2q}}^2\bigg)^{\frac12}
  +\sum_{i=1}^m\|\nabla A_i\|_{L^{2q}}^2\bigg]\\
  &\hskip26mm+\frac{1}{\sigma^2}\sum_{i=1}^m\|A_i-\hat
  A_i\|_{L^{2q}}^2
  +\frac{1}\sigma\bigg[\|A_0-\hat A_0\|_{L^q}+\bigg(\sum_{i=1}^m\|A_i-\hat
  A_i\|_{L^{2q}}^2\bigg)^{\frac12}\bigg]\bigg\},
  \end{split}
  \end{equation*}
where $p$ is the conjugate number of $q$: $1/p+1/q=1$ and
\begin{equation}\label{5.4}
G_R(w)=\bigg\{x\in\R^d:\sup_{0\leq t\leq T}|X_t(w,x)|\vee|\hat
X_t(w,x)|\leq R\bigg\}.
\end{equation}
\end{theorem}

\noindent{\bf Proof.} Denote by $\xi_t=X_t-\hat X_t$, then
$\xi_0=0$. By It\^{o} formula,
  \begin{eqnarray}\label{5.5}
  \d|\xi_t|^2&=&2\sum_{i=1}^m\langle\xi_t,A_i(X_t)-\hat A_i(\hat
  X_t)\rangle\d w^i_t+2\langle\xi_t,A_0(X_t)-\hat A_0(\hat
  X_t)\rangle\d t\cr
  &&+\sum_{i=1}^m|A_i(X_t)-\hat A_i(\hat X_t)|^2\d t.
  \end{eqnarray}
For $\sigma>0$, $\log\big(\frac{|\xi_t|^2}{\sigma^2}+1\big)=
\log(|\xi_t|^2+\sigma^2)-\log \sigma^2$. Again by the It\^{o}
formula,
  $$\d\log(|\xi_t|^2+\sigma^2)=\frac{\d|\xi_t|^2}{|\xi_t|^2+\sigma^2}
  -\frac12\cdot\frac{4\sum_{i=1}^m\langle\xi_t,A_i(X_t)-\hat A_i(\hat
  X_t)\rangle^2}{(|\xi_t|^2+\sigma^2)^2}\d t,$$
using \eqref{5.5}, we obtain
  \begin{eqnarray}\label{5.6}
  \d\log(|\xi_t|^2+\sigma^2)&=&2\sum_{i=1}^m\frac{\langle\xi_t,A_i(X_t)-\hat A_i(\hat
  X_t)\rangle}{|\xi_t|^2+\sigma^2}\d w^i_t
  +2\frac{\langle\xi_t,A_0(X_t)-\hat A_0(\hat
  X_t)\rangle}{|\xi_t|^2+\sigma^2}\d t\cr
  &&+\sum_{i=1}^m\frac{|A_i(X_t)-\hat A_i(\hat
  X_t)|^2}{|\xi_t|^2+\sigma^2}\d t
  -2\sum_{i=1}^m\frac{\langle\xi_t,A_i(X_t)-\hat A_i(\hat
  X_t)\rangle^2}{(|\xi_t|^2+\sigma^2)^2}\d t\cr
  &=:&\d I_1(t)+\d I_2(t)+\d I_3(t)+\d I_4(t).
  \end{eqnarray}

Let $\tau_R(x)=\inf\{t\geq0:|X_t(x)|\vee|\hat X_t(x)|> R\}$. Remark
that almost surely, $G_R\subset\{x:\tau_R(x)>T\}$ and for any
$t\geq0$, $\{\tau_R>t\}\subset B(R)$. Therefore
  \begin{equation*}
  \E\bigg[\int_{G_R}\sup_{0\leq t\leq T}|I_1(t)|\,\d\gamma_d\bigg]
  \leq  \E\bigg[\int_{B(R)}\sup_{0\leq t\leq T\wedge\tau_R}|I_1(t)|\,\d\gamma_d\bigg].
  \end{equation*}
By Burkholder's inequality,
  \begin{equation*}
  \E\bigg(\sup_{0\leq t\leq T\wedge\tau_R}|I_1(t)|^2\bigg)\leq
  4\,\E\biggl(\int_0^{T\wedge\tau_R}\sum_{i=1}^m\frac{\langle\xi_t,A_i(X_t)-\hat
  A_i(\hat X_t)\rangle^2}{(|\xi_t|^2+\sigma^2)^2}\,\d t\bigg),
  \end{equation*}
which is obviously less than
  \begin{equation*}
  4\,\E\biggl(\int_0^{T\wedge\tau_R}\sum_{i=1}^m\frac{|A_i(X_t)-\hat
  A_i(\hat X_t)|^2}{|\xi_t|^2+\sigma^2}\,\d t\bigg).
  \end{equation*}
Hence
  \begin{eqnarray}\label{5.7}
  \hskip-3mm\E\bigg[\int_{B(R)}\sup_{0\leq t\leq
  T\wedge\tau_R}|I_1(t)|\,\d\gamma_d\bigg]\leq
  4\,\bigg[\int_0^{T}\bigg(\E\int_{\{\tau_R>t\}}\sum_{i=1}^m\frac{|A_i(X_t)-\hat A_i(\hat
  X_t)|^2}{|\xi_t|^2+\sigma^2}\d\gamma_d\,\bigg)\d t\bigg]^{\frac12}.
  \end{eqnarray}
We have $A_i(X_t)-\hat A_i(\hat X_t)=A_i(X_t)-A_i(\hat X_t)+A_i(\hat
X_t)-\hat A_i(\hat X_t)$. Using the density $\hat K_t$, it is clear
that
  \begin{eqnarray*}
  \E\int_{\{\tau_R>t\}}\frac{|A_i(\hat X_t)-\hat A_i(\hat
  X_t)|^2}{|\xi_t|^2+\sigma^2}\d\gamma_d
  &\leq&\frac1{\sigma^2}\E\int_{\R^d}|A_i(\hat X_t)-\hat A_i(\hat
  X_t)|^2\d\gamma_d\cr
  &=&\frac1{\sigma^2}\E\int_{\R^d}|A_i-\hat A_i|^2\hat
  K_t\,\d\gamma_d.
  \end{eqnarray*}
Thus by H\"{o}lder's inequality and according to \eqref{2.3}, we
have
  \begin{eqnarray}\label{5.8}
  \E\int_{\{\tau_R>t\}}\frac{|A_i(\hat X_t)-\hat A_i(\hat
  X_t)|^2}{|\xi_t|^2+\sigma^2}\d\gamma_d
  \leq\frac{\Lambda_{p,T}}{\sigma^2}\|A_i-\hat A_i\|_{L^{2q}}^2.
  \end{eqnarray}

Now we shall use Theorem \ref{th.6.2} in the Appendix to estimate
another term. Note that on the set $\{\tau_R>t\}$, $X_t,\hat X_t\in
B(R)$, then $|X_t-\hat X_t|\leq 2R$. Since
$(X_t)_\#\gamma_d\ll\gamma_d$ and $(\hat
X_t)_\#\gamma_d\ll\gamma_d$, we can apply \eqref{th.6.2-1} so that
  \begin{equation*}
  |A_i(X_t)-A_i(\hat X_t)|\leq C_d |X_t-\hat X_t|\,\bigl(
  M_{2R}|\nabla A_i|(X_t)+M_{2R}|\nabla A_i|(\hat X_t)\bigr).
  \end{equation*}
Then
  \begin{equation*}
  \E\bigg[\int_{\{\tau_R>t\}}\frac{|A_i(X_t)-A_i(\hat
  X_t)|^2}{|\xi_t|^2+\sigma^2}\d\gamma_d\bigg]
  \leq  C_d^2\,\E\int_{\{\tau_R>t\}}\big(M_{2R}|\nabla A_i|(X_t)
  +M_{2R}|\nabla A_i|(\hat
  X_t)\big)^2\,\d\gamma_d.
  \end{equation*}
Notice again that on $\{\tau_R(x)>t\}$,  $X_t(x)$ and $\hat X_t(x)$
are in $B(R)$, therefore
  \begin{eqnarray}\label{5.8.5}
  \E\bigg[\int_{\{\tau_R>t\}}\frac{|A_i(X_t)-A_i(\hat
  X_t)|^2}{|\xi_t|^2+\sigma^2}\d\gamma_d\bigg]
  &\leq&2C_d^2\,\E\int_{B(R)}(M_{2R}|\nabla A_i|)^2\big(K_t+\hat
  K_t\big)\d\gamma_d\cr
  &\leq& 4C_{d}^2\Lambda_{p,T}\bigg(\int_{B(R)}\big(M_{2R}|\nabla A_i|\big)^{2q}\d\gamma_d\bigg)^{\frac1q}.
  \end{eqnarray}

Remark that the maximal function inequality does not hold for the
Gaussian measure $\gamma_d$ on the whole space $\R^d$. However, on
each ball $B(R)$,
\begin{equation*}
\gamma_d|_{B(R)}\leq \frac{1}{(2\pi)^{d/2}}\Leb_d|_{B(R)}\leq
e^{R^2/2}\gamma_d|_{B(R)}.
\end{equation*}
Thus, according to \eqref{th.6.2-2},
  \begin{equation*}
  \begin{split}
  \int_{B(R)}\big(M_{2R}|\nabla A_i|\big)^{2q}\d\gamma_d
  &\leq \frac1{(2\pi)^{d/2}}\int_{B(R)}\big(M_{2R}|\nabla
  A_i|\big)^{2q}\d x
  \leq\frac{C_{d,q}}{(2\pi)^{d/2}}\int_{B(3R)}|\nabla
  A_i|^{2q}\d x\\
  &\leq C_{d,q}e^{9R^2/2}\int_{B(3R)}|\nabla
  A_i|^{2q}\d\gamma_d
  \leq C_{d,q}e^{9R^2/2}\|\nabla A_i\|_{L^{2q}}^{2q}.
  \end{split}
  \end{equation*}
Therefore by \eqref{5.8.5}, there exists a constant $C_{d,q,R}>0$
such that
  \begin{eqnarray*}
  \E\bigg[\int_{\{\tau_R>t\}}\frac{|A_i(X_t)-A_i(\hat
  X_t)|^2}{|\xi_t|^2+\sigma^2}\d\gamma_d\bigg]
  \leq C_{d,q,R}\Lambda_{p,T}\|\nabla A_i\|_{L^{2q}}^{2}.
  \end{eqnarray*}
Combining this estimate with \eqref{5.7} and \eqref{5.8}, we get
  \begin{equation}\label{5.10}
  \E\bigg[\int_{G_R}\sup_{0\leq t\leq T}|I_1(t)|\d\gamma_d\bigg]
  \leq CT^{\frac12}\Lambda_{p,T}^{\frac12}\bigg(C_{d,q,R}\sum_{i=1}^m\|\nabla A_i\|_{L^{2q}}^2
  +\frac1{\sigma^2}
  \sum_{i=1}^m\|A_i-\hat A_i\|_{L^{2q}}^2\bigg)^{\frac12}.
  \end{equation}

Now we turn to deal with $I_2(t)$ in \eqref{5.6}. We have
  \begin{equation*}
  \E\bigg[\int_{G_R}\sup_{0\leq t\leq
  T}|I_2(t)|\,\d\gamma_d\bigg]
  \leq 2\int_0^T\bigg[\E\int_{G_R}\frac{|A_0(X_t)-\hat A_0(\hat X_t)|}
  {(|\xi_t|^2+\sigma^2)^{\frac12}}\d\gamma_d\bigg]\,\d t.
  \end{equation*}
Note that for $x\in G_R$, $\hat X_t(x)\in B(R)$ for each
$t\in[0,T]$, thus
  \begin{equation*}
  \E\bigg[\int_{G_R}\frac{|A_0(\hat X_t)-\hat A_0(\hat X_t)|}
  {(|\xi_t|^2+\sigma^2)^{\frac12}}\d\gamma_d\bigg]
  \leq \frac1\sigma\E\int_{B(R)}|A_0-\hat A_0|\,\hat K_t\,\d\gamma_d
  \leq \frac{\Lambda_{p,T}}{\sigma}\,\|A_0-\hat A_0\|_{L^q}.
  \end{equation*}
Again using \eqref{th.6.2-1},
  \begin{equation*}
  \E\bigg[\int_{G_R}\frac{|A_0(X_t)-A_0(\hat X_t)|}
  {(|\xi_t|^2+\sigma^2)^{\frac12}}\d\gamma_d\bigg]
  \leq C_d\,\E\int_{G_R}\big(M_{2R}|\nabla A_0|(X_t)+M_{2R}|\nabla A_0|(\hat
  X_t)\big)\d\gamma_d,
  \end{equation*}
which is dominated by
  \begin{equation*}
  C_d\, \E\bigg[\int_{B(R)} (M_{2R}|\nabla
  A_0|)\cdot(K_t+ \hat K_t)\,d\gamma_d\bigg] \leq C_{d,q,R}\,\|\nabla
  A_0\|_{L^q}\,\Lambda_{p,T}.
  \end{equation*}
Therefore we get the following estimate for $I_2$:
  \begin{eqnarray}\label{5.11}
  \E\bigg[\int_{G_R}\sup_{0\leq t\leq
  T}|I_2(t)|\,\d\gamma_d\bigg]
  \leq 2T\Lambda_{p,T}\bigg(C_{d,q,R}\|\nabla
  A_0\|_{L^q}+\frac1\sigma\|A_0-\hat
  A_0\|_{L^q}\bigg).
  \end{eqnarray}

In the same way we have
  \begin{equation}\label{5.12}
  \E\bigg[\int_{G_R}\sup_{0\leq t\leq T}|I_3(t)|\,\d\gamma_d\bigg]
  \leq CT\Lambda_{p,T}\bigg(C_{d,q,R}\sum_{i=1}^m\|\nabla A_i\|_{L^{2q}}^2
  +\frac1{\sigma^2}
  \sum_{i=1}^m\|A_i-\hat A_i\|_{L^{2q}}^2\bigg).
  \end{equation}
The term $I_4(t)$ is negative and hence we omit it. Combining
\eqref{5.6} and \eqref{5.10}--\eqref{5.12}, we complete the proof.
\fin

\medskip

Now we shall construct a solution to SDE \eqref{5.1}. To this end,
we take $\ee=1/n$ and we write $A_j^n$ instead of $A_j^{1/n}$
introduced in Section 3. Then by assumption (A2) and Lemma
\ref{sect-3-lem-1}, there is $C>0$ independent of $n$ and $i$, such
that
  \begin{equation}\label{5.13}
  |A^n_i(x)|\leq C(1+|x|).
  \end{equation}
Let $X^n_t$ be the solution to It\^{o} SDE \eqref{5.1} with the
coefficients $A^n_i\,(i=0,1,\ldots, m)$. Then for any $\alpha\geq1$
and $T>0$, there exists $C_{\alpha,T}>0$ independent of $n$ such
that
  \begin{equation}\label{5.14}
  \E\bigg(\sup_{0\leq t\leq T}|X^n_t|^\alpha\bigg)\leq
  C_{\alpha,T}(1+|x|^\alpha),\quad\mbox{for all }x\in\R^d.
  \end{equation}
Let $K^n_t$ be the density of $(X^n_t)_\#\gamma_d$ with respect to
$\gamma_d$. Under the hypotheses (A2)--(A4), there is $T_0>0$ such
that (recall that $p$ is the conjugate number of $q>1$):
  \begin{eqnarray}\label{LambdapT0}
  \hskip-8mm\Lambda_{p,T_0}&:=&\bigg[\int_{\R^d}\exp\bigg(2pT_0\Big[|A_0|+e|\delta(A_0)|\cr
  &&\hskip18mm+\sum_{j=1}^m\big(2p|A_j|^2+|\nabla
  A_j|^2+2(p-1)e^2|\delta(A_j)|^2\big)\Big]\bigg)\d\gamma_d\bigg]^{\frac{p-1}{p(2p-1)}}<\infty.
  \end{eqnarray}
Similar to \eqref{LambdaT0New}, we have
  \begin{equation}\label{5.15}
  \sup_{t\in[0,T_0]}\sup_{n\geq1}\|K^n_t\|_{L^p(\gamma_d\times\P)}\leq \Lambda_{p,T_0}<+\infty.
  \end{equation}

Now we shall prove that the family $\{X^n_\cdot:n\geq1\}$ is
convergent to some stochastic field.

\begin{theorem}\label{th.5.2}
Let $T_0$ be given in \eqref{LambdapT0}. Then under the assumptions
\textup{(A1')} and \textup{(A2)--(A4)}, there exists
$X:\Omega\times\R^d\ra C([0,T_0],\R^d)$ such that for any
$\alpha\geq 1$,
\begin{equation}\label{5.16}
  \lim_{n\ra\infty}\E\bigg[\int_{\R^d}\bigg(\sup_{0\leq t\leq
  T_0}|X^n_t-X_t|^\alpha\bigg)\d\gamma_d\bigg]=0.
\end{equation}
\end{theorem}

\noindent{\bf Proof.} We shall prove that $\{X^n;\ n\geq 1\}$ is a
Cauchy sequence in $L^\alpha(\Omega\times\R^d; C([0,T_0],\R^d))$.
Denote by $\|\cdot\|_{\infty,T_0}$ the uniform norm on
$C([0,T_0],\R^d)$, so what we have to prove is
  \begin{equation}\label{th.5.2.1}
  \lim_{n,k\ra +\infty}\E\bigg(\int_{\R^d}
  \|X^n-X^k\|_{\infty,T_0}^\alpha\,\d\gamma_d\bigg)=0.
  \end{equation}

First by \eqref{5.14}, the quantity
  \begin{equation}\label{th.5.2.2}
  J_{\alpha,T_0}:=\sup_{n\geq
  1}\E\bigg(\int_{\R^d}\|X^n\|_{\infty,T_0}^{2\alpha}\,\d\gamma_d\bigg)\leq
  C_{\alpha,T_0}\,\int_{\R^d} (1+|x|^{2\alpha})\,\d\gamma_d
  \end{equation}
is obviously finite. Let $R>0$ and set
  \begin{equation*}
  G_{n,R}(w)=\{x\in\R^d;\ \|X^n(w,x)\|_{\infty,T_0}\leq
  R\}.
  \end{equation*}
Using \eqref{th.5.2.2}, for any $\alpha\geq 1$ and $R>0$, we have
  \begin{equation*}
  \sup_{n\geq1}\E\bigl(\gamma_d(G_{n,R}^c)\bigr)\leq \frac{J_{\alpha,T_0}}{
  R^{2\alpha}}.
  \end{equation*}
Now by Cauchy-Schwarz inequality
  \begin{equation*}
  \begin{split}
  &\E\biggl(\int_{G_{n,R}^c\cup
  G_{k,R}^c}\|X^n-X^k\|_{\infty,T_0}^\alpha\,\d\gamma_d\biggr)\\
  &\hskip6mm\leq \Big(\E\bigl[\gamma_d\big(G_{n,R}^c\cup
  G_{k,R}^c\big)\bigr]\Big)^{1/2}\cdot
  \biggl(\E\int_{\R^d}\|X^n-X^k\|_{\infty,T_0}^{2\alpha}\,\d\gamma_d\biggr)^{1/2}\\
  &\hskip6mm\leq \bigg(\frac{2J_{\alpha,T_0}}{ R^{2\alpha}}\bigg)^{1/2}\cdot
  \bigl(2^{2\alpha} J_{\alpha,T_0}\bigr)^{1/2}.
  \end{split}
  \end{equation*}
Let $\ee>0$ be given; choose $R>1$ big enough such that the last
quantity in the above inequality is less than $\ee$. Then we have
for any $n,k\geq 1$,
  \begin{equation}\label{th.5.2.3}
  \E\biggl(\int_{G_{n,R}^c\cup
  G_{k,R}^c}\|X^n-X^k\|_{\infty,T_0}^\alpha\,\d\gamma_d\biggr)\leq\ee.
  \end{equation}

Let
  \begin{equation*}
  \sigma_{n,k}=\|A_0^n-A_0^k\|_{L^q}+\bigg(\sum_{i=1}^m\|A_i^n-A_i^k\|_{L^{2q}}^2\bigg)^{1/2},
  \end{equation*}
which tends to 0 as $n,k\ra +\infty$ since $A_0^n$ converges to
$A_0$ in $L^q(\gamma_d)$ and $A_i^n$ converges to $A_i$ in
$L^{2q}(\gamma_d)$ for $i=1, \cdots, m$. Now applying Theorem
\ref{th.5.1} with $A_i$ and  $\hat A_i$ being replaced respectively
by $A_i^n$ and $A_i^k$, we get
  \begin{eqnarray*}
  I_{n,k}&:=&\E\biggl[\int_{G_{n,R}\cap
  G_{k,R}}\log\bigg(\frac{\|X^n-X^k\|_{\infty,T_0}^2}{\sigma_{n,k}^2}+1\bigg)\,\d\gamma_d\biggr]\cr
  &\leq& C_{T_0}\Lambda_{p,T_0}\bigg\{C_{d,q,R}\bigg[\|\nabla
  A_0^n\|_{L^q}+ \bigg(\sum_{i=1}^m\|\nabla
  A_i^n\|_{L^{2q}}^2\bigg)^{1/2}+\sum_{i=1}^n\|\nabla
  A_i^n\|_{L^{2q}}^2\bigg]+2\bigg\}.
  \end{eqnarray*}
Recall that $A_i^n=\varphi_{1/n}\,P_{1/n}A_i$, then $\nabla A_i^n=
\nabla\varphi_{1/n}\otimes P_{1/n}A_i+ \varphi_{1/n}\,
e^{-1/n}P_{1/n} \nabla A_i$, therefore
  $$|\nabla A_i^n|\leq P_{1/n}(|A_i|+ |\nabla A_i|).$$
We get the following uniform estimates
  \begin{equation*}
  \|\nabla A_0^n\|_{L^{q}}\leq \|A_0\|_{\D_1^q},\quad \|\nabla
  A_i^n\|_{L^{2q}}\leq \|A_i\|_{\D_1^{2q}}.
  \end{equation*}
So the quantity $I_{n,k}$ is uniformly bounded with respect to
$n,k$. Let $\hat\Pi$ be the measure on $\Omega\times\R^d$ defined by
  \begin{equation*}
  \int_{\Omega\times\R^d}
  \psi(w,x)\,\d\hat\Pi(w,x)=\E\bigg[\int_{G_{n,R}\cap
  G_{k,R}}\psi(w,x)\,\d\gamma_d(x)\bigg].
  \end{equation*}
We have $\hat\Pi(\Omega\times\R^d)\leq 1$. Let $\eta>0$, consider
  \begin{equation*}
  \Sigma_{n,k}=\{(w,x);\ \|X^n(w,x)-X^k(w,x)\|_{\infty,T_0}\geq\eta\},
  \end{equation*}
which is equal to
  \begin{equation*}
  \biggl\{(w,x);\
  \log\bigg(\frac{\|X^n-X^k\|_{\infty,T_0}^2}{\sigma_{n,k}^2}+1\bigg)\geq
  \log\bigg(\frac{\eta^2}{\sigma_{n,k}^2}+1\bigg)\biggr\}.
  \end{equation*}
It follows that as $n,k\ra+\infty$,
  \begin{equation}\label{A2.1}
  \hat\Pi(\Sigma_{n,k})\leq \frac{I_{n,k}}{
  \log\Bigl(\frac{\eta^2}{\sigma_{n,k}^2}+1\Bigr)}\ra 0,
  \end{equation}
since $\sigma_{n,k}\ra0$ and the family $\{I_{n,k};n,k\geq 1\}$ is
bounded. Now
  \begin{eqnarray}\label{th.5.2.5}
  &&\E\bigg(\int_{G_{n,R}\cap
  G_{k,R}}\|X^n-X^k\|_{\infty,T_0}^\alpha\,\d\gamma_d\bigg)
  =\int_{\Omega\times\R^d}\|X^n-X^k\|_{\infty,T_0}^\alpha\,\d\hat\Pi\cr
  &&\hskip6mm=\int_{\Sigma_{n,k}^c}\|X^n-X^k\|_{\infty,T_0}^\alpha\,\d\hat\Pi
  +\int_{\Sigma_{n,k}}\|X^n-X^k\|_{\infty,T_0}^\alpha\,\d\hat\Pi.
  \end{eqnarray}
The first term on the right side of \eqref{th.5.2.5} is less than
$\eta^\alpha$, while the second one, due to \eqref{th.5.2.2} and
\eqref{A2.1}, is dominated by
  \begin{eqnarray*}
  \sqrt{\hat\Pi(\Sigma_{n,k})}\cdot
  \sqrt{\E\int_{\R^d}\|X^n-X^k\|_{\infty,T_0}^{2\alpha}\d\gamma_d}
  \leq 2^{\alpha}\sqrt{J_{\alpha,T_0}\hat\Pi(\Sigma_{n,k})}\ra 0\quad
  \hbox{as }n,k\ra+\infty.
  \end{eqnarray*}
Now taking $\eta=\ee^{1/\alpha}$ and combining \eqref{th.5.2.3} and
\eqref{th.5.2.5}, we prove that
  \begin{equation*}
  \limsup_{n,k\ra
  +\infty}\E\bigg[\int_{\R^d}\|X^n-X^k\|_{\infty,T_0}^\alpha\,\d\gamma_d\bigg]\leq
  2\ee,
  \end{equation*}
which implies \eqref{th.5.2.1}.

Let $X\in L^\alpha(\Omega\times\R^d; C([0,T_0],\R^d))$ be the limit
of $X^n$ in this space. We see that for each $t\in [0,T]$ and almost
all $x\in\R^d$, $w\ra X_t(w,x)$ is in $\F_t$.  \fin

\begin{proposition}\label{prop5.2}
There exists a family $\{\hat K_t; t\in[0,T_0]\}$ of density
functions on $\R^d$ such that $(X_t)_\#\gamma_d=\hat K_t \gamma_d$
for each $t\in[0,T_0]$. Moreover, $\sup_{0\leq t\leq T_0}\|\hat
K_t\|_{L^p(\P\times\gamma_d)}\leq \Lambda_{p,T_0}$, where
$\Lambda_{p,T_0}$ is given in \eqref{5.15}.
\end{proposition}

\noindent{\bf Proof.} It is the same as the proof of Theorem
\ref{sect-3-thm-1}. \fin

\medskip

The same arguments in the proof of Proposition \ref{prop4.1} and
\ref{prop4.2} yield the following

\begin{proposition}\label{prop5.3} For any $\alpha\geq 2$, up to a
subsequence,
  \begin{equation*}
  \lim_{n\ra\infty}\int_{\R^d}\E\bigg(\sup_{0\leq t\leq T_0}
  \bigg|\sum_{i=1}^m\int_0^t\big[A^n_i(X^n_s)-A_i(X_s)\big]\d
  w^i_s\bigg|^\alpha\bigg)\d\gamma_d=0,
  \end{equation*}
and
  \begin{equation*}
  \lim_{n\ra\infty}\int_{\R^d}\biggl[\E\int_0^{T_0}
  |A^n_0(X^n_s)-A_0(X_s)|^\alpha\d s\biggr]\d\gamma_d=0.
  \end{equation*}
\end{proposition}

Now for regularized vector fields $A_i^n, i=0,1, \cdots, m$, we have
  \begin{equation}\label{5.23}
  X^n_t(x)=x+\sum_{i=1}^m\int_0^tA^n_i(X^n_s)\,\d w^i_s
  +\int_0^tA^n_0(X^n_s)\,\d s.
  \end{equation}
When $n\ra +\infty$, by Theorem \ref{th.5.2} and Proposition
\ref{prop5.3}, the two sides of \eqref{5.23} converge respectively
to $X$ and
  \begin{equation*}
  x+\sum_{i=1}^m\int_0^\cdot A_i(X_s)\,\d w^i_s
  +\int_0^\cdot A_0(X_s)\,\d s
  \end{equation*}
in the space $L^\alpha(\Omega\times\R^d; C([0,T_0],\R^d))$.
Therefore for almost all $x\in\R^d$, the following equality holds
$\P$-almost surely:
  $$X_t(x)=x+\sum_{i=1}^m\int_0^tA_i(X_s)\,\d w^i_s
  +\int_0^tA_0(X_s)\,\d s,\quad\mbox{for all }t\in[0,T_0].$$
That is to say, $X_t$ solves SDE \eqref{5.1} over $[0,T_0]$.

The following result proves the pathwise uniqueness to SDE
\eqref{5.1} for a.e. initial value $x\in\R^d$.

\begin{proposition}\label{prop5.4}
Under the conditions \textup{(A1')} and \textup{(A2)--(A4)}, the SDE
\eqref{5.1} has a unique solution on the interval $[0,T_0]$.
\end{proposition}

\noindent{\bf Proof.} Let $(Y_t)_{t\in[0,T_0]}$ be another solution.
Set, for $R>0$,
  \begin{equation*}
  G_R=\bigg\{(w,x)\in \Omega\times\R^d;\ \sup_{0\leq t\leq
  T_0}|X_t(w,x)-Y_t(w,x)|\leq R\bigg\}.
  \end{equation*}
Remark that in Theorem \ref{th.5.1}, the terms involving $1/\sigma$
and $1/\sigma^2$ are equal to zero. Therefore the term
  \begin{eqnarray*}
  I&:=&\E\int_{G_R}\log\bigg(\frac{\sup_{0\leq t\leq
  T_0}|X_t-Y_t|^2}{\sigma^2}+1\bigg)\d\gamma_d\cr
  &\leq& C_{T_0}\Lambda_{p,T_0}C_{d,q,R}
  \bigg[\|A_0\|_{\D_1^q}+\bigg(\sum_{i=1}^m\|A_i\|_{\D_1^{2q}}^2\bigg)^{\frac12}
  +\sum_{i=1}^m\|A_i\|_{\D_1^{2q}}^2\bigg]
  \end{eqnarray*}
is bounded for any $\sigma>0$. Consider for $\eta>0$,
  \begin{equation*}
  \Sigma_\eta=\bigg\{(w,x); \sup_{0\leq t\leq T_0}|X_t(w,x)-
  Y_t(w,x)|\geq\eta\bigg\}.
  \end{equation*}
Similar to \eqref{A2.1}, we have
  \begin{equation*}
  \E\bigg(\int_{G_R} {\bf 1}_{\Sigma_\eta}\d\gamma_d\bigg) \leq
  \frac{I}{\log\bigl(\frac{\eta^2}{\sigma^2}+1\bigr)}\ra 0
  \end{equation*}
as $\sigma\ra 0$. So we obtain
  $$\ch_{G_R}\cdot\sup_{0\leq t\leq T_0}|X_t-Y_t|=0,\quad (\P\times \gamma_d)\mbox{-a.s.}$$
Letting $R\ra\infty$, we obtain that $(\P\times\gamma_d)$ almost
surely, $X_t=Y_t$ for all $t\in[0,T_0]$. \fin

\medskip

Now we extend the solution to any time interval $[0,T]$. Let
$\theta_{T_0}w$ be the time-shift of the Brownian motion $w$ and
denote by $X_t^{T_0}$ the corresponding solution to SDE driven by
$\theta_{T_0}w$. By Proposition \ref{prop5.4},
$\{X_t^{T_0}(\theta_{T_0}w,x): 0\leq t\leq T_0\}$ is the unique
solution to the SDE over $[0,T_0]$:
  \begin{equation*}
  X_t^{T_0}(x)=x+\sum_{i=1}^m \int_0^t A_i(X_s^{T_0}(x))\,\d(\theta_{T_0}w)^i_s
  +\int_0^tA_0(X_s^{T_0}(x))\,\d s.
  \end{equation*}

For $t\in[0,T_0]$, define $X_{t+T_0}(w,x)=X^{T_0}_t(\theta_{T_0}w,
X_{T_0}(w,x))$. Note that $X_t$ is well defined on the interval
$[0,2T_0]$ up to a $(\P\times\gamma_d)$-negligible subset of
$\Omega\times\R_d$. Replacing $x$ by $X_{T_0}(x)$ in the above
equation, we get easily
  \begin{equation*}
  X_{t+T_0}(x)=x+\sum_{i=1}^m \int_0^{t+T_0} A_i(X_s(x))\,\d w^i_s
  +\int_0^{t+T_0}A_0(X_s(x))\,\d s.
  \end{equation*}
Therefore $X_t$ defined as above is a solution to SDE on the
interval $[0,2T_0]$. Continuing in this way, we obtain the solution
of SDE \eqref{5.1} on $[0,T]$.

\begin{theorem}\label{th.5.3}
The $\{X_t;\ t\in [0,T]\}$ constructed above is the unique solution
to SDE \eqref{5.1} in the sense of Definition \ref{def.5.1}.
Moreover for each $t\in [0,T]$, the density $K_t$ of
$(X_t)_\#\gamma_d$ with respect to $\gamma_d$ is in the space
$L^1\log{L^1}$.
\end{theorem}

\noindent{\bf Proof.} Let $Y_t,\,t\in[0,T]$ be another solution in
the sense of Definition \ref{def.5.1}. First by Proposition
\ref{prop5.4}, we have $(\P\times \gamma_d)$-almost surely,
$Y_t=X_t$ for all $t\in [0,T_0]$. In particular, $Y_{T_0}=X_{T_0}$.
Next by the flow property, $Y_{t+T_0}$ satisfies the following
equation:
  \begin{equation*}
  Y_{t+T_0}(x)=Y_{T_0}(x)+\sum_{i=1}^m \int_0^t A_i\big(Y_{s+T_0}(x)\big)\d(\theta_{T_0}w)^i_s
  +\int_0^tA_0\big(Y_{s+T_0}(x)\big)\d s,
  \end{equation*}
that is, $Y_{t+T_0}$ is a solution with initial value $Y_{T_0}$. But
by the above discussion, $X_{t+T_0}$ is also a solution with the
same initial value $X_{T_0}=Y_{T_0}$. Again by Proposition
\ref{prop5.4}, we have $(\P\times\gamma_d)$-almost surely,
$X_{t+T_0}=Y_{t+T_0}$ for all $t\leq T_0$. Hence we have proved
$X|_{[0,2T_0]}=Y|_{[0,2T_0]}$. Repeating this procedure, we obtain
the uniqueness over $[0,T]$. The existence of density $K_t$ of
$(X_t)_\#\gamma_d$ with respect to $\gamma_d$ beyond $T_0$ is
deduced from the flow property. However, to insure that $K_t\in
L^1\log{L^1}$, we have to use Theorem \ref{sect-3-prop-1} and the
following
  \begin{equation*}
  \lim_{n\ra\infty}\int_{\R^d}\E\bigg(\sup_{0\leq t\leq T}
  |X^n_t-X_t|^\alpha\bigg)\d\gamma_d=0,
  \end{equation*}
which can be checked using the same arguments as in the proof of
Propositions \ref{prop4.1} and \ref{prop4.2}. \fin

\section{Appendix}

For any locally integrable function $f\in L_{loc}^1(\R^d)$ and
$R>0$, the local maximal function $M_Rf$ is defined by
  \begin{equation}\label{6.1}
  M_R f(x)=\sup_{0<r\leq R}\frac{1}{\Leb_d(B(x,r))}\int_{B(x,r)}|f(y)|\,\d y,
  \end{equation}
where $B(x,r)=\{y\in\R^d; |y-x|\leq r\}$. The following result is
the starting point for the approach concerning Sobolev coefficients,
used in \cite{CrippadeLellis} and \cite{Zhang09}.

\begin{theorem}\label{th.6.2}
Let $f\in L_{loc}^1(\R^d)$ be such that $\nabla f\in
L_{loc}^1(\R^d)$. Then there is a constant $C_d>0$ (independent of
$f$) and a negligible subset $N$, such that for $x,y\in N^c$ with
$|x-y|\leq R$,
  \begin{equation}\label{th.6.2-1}
  |f(x)-f(y)|\leq C_d |x-y|\, \bigl((M_R|\nabla f|)(x)+(M_R|\nabla
  f|)(y)\bigr);
  \end{equation}
moreover for $p>1$ and $f\in L_{loc}^p(\R^d)$, there is a constant
$C_{d,p}>0$ such that
  \begin{equation}\label{th.6.2-2}
  \int_{B(r)} (M_R f)^p\,\d x\leq C_{d,p}\, \int_{B(r+R)}|f|^p\,\d x.
  \end{equation}
\end{theorem}

Since the inequality \eqref{th.6.2-1} played a key role in the proof
of Theorem \ref{th.5.1}, we give here its proof for the sake of the
reader's convenience.

We follow the idea of the proof of Claim \#2 on p.253 in
\cite{EvansGariepy}. For any bounded measurable subset $U$ in $\R^d$
such that its Lebesgue measure $\Leb_d(U)>0$, define the average of
$f\in L^1_{loc}(\R^d)$ on $U$ by
  $$(f)_U=\bint_{U}f(y)\,\d y
  :=\frac1{\Leb_d(U)}\int_{U} f(y)\,\d y.$$
Write $(f)_{x,r}$ instead of $(f)_{B(x,r)}$ for simplicity. Then
$M_Rf(x)=\sup_{0<r\leq R}(|f|)_{x,r}$. We will need the following
simple inequality: for any $C\in\R$,
  \begin{equation}\label{I0}
  |(f)_{U}-C|\leq\bint_{U}|f(y)-C|\,\d y.
  \end{equation}

First, for any $x\in\R^d$ and $r\in\,]0,R]$, by Poincar\'{e}'s
inequality with $p=1$ and $p^\ast=d/(d-1)$ (see \cite{EvansGariepy}
p.141), there is $C_d>0$ such that
  \begin{eqnarray}\label{II}
  \bint_{B(x,r)}|f-(f)_{x,r}|\,\d y
  &\leq&\bigg(\bint_{B(x,r)}|f-(f)_{x,r}|^{d/(d-1)}\,\d y\bigg)^{(d-1)/d}\cr
  &\leq& C_d\,r\bint_{B(x,r)}|\nabla f|\,\d y
  \leq C_dM_R|\nabla f|(x)\,r.
  \end{eqnarray}
In particular, for all $k\geq 0$, by \eqref{I0} and \eqref{II},
  \begin{eqnarray*}
  |(f)_{x,r/2^{k+1}}-(f)_{x,r/2^k}| &\leq&
  \bint_{B(x,r/2^{k+1})}|f-(f)_{x,r/2^k}|\,\d y\cr
  &\leq&2^d\bint_{B(x,r/2^k)}|f-(f)_{x,r/2^k}|\,\d y\cr
  &\leq&2^dC_dM_R|\nabla f|(x)\,r/2^k.
  \end{eqnarray*}
Since $f\in L^1_{loc}(\R^d)$, there is a negligible subset $N\subset
\R^d$, such that for all $x\in N^c$, $f(x)=\lim_{r\ra0}(f)_{x,r}$.
Thus for any $x\in N^c$, we have by summing up the above inequality
that
  \begin{eqnarray}\label{III}
  |f(x)-(f)_{x,r}|\leq
  \sum_{k=0}^\infty|(f)_{x,r/2^{k+1}}-(f)_{x,r/2^k}|
  \leq2^{1+d}C_dM_R|\nabla f|(x)\,r.
  \end{eqnarray}

Next for all $x,y\in N^c, x\neq y$ and $|x-y|\leq R$, let $r=|x-y|$.
Then by the triangular inequality, \eqref{I0} and \eqref{II},
  \begin{eqnarray}\label{IV}
  |(f)_{x,r}-(f)_{y,r}|&\leq&
  \bint_{B(x,r)\cap
  B(y,r)}\big(|(f)_{x,r}-f(z)|+|f(z)-(f)_{y,r}|\big)\d z\cr
  &\leq &\tilde C_d\bigg[\bint_{B(x,r)}|(f)_{x,r}-f(z)|\,\d z
  +\bint_{B(y,r)}|f(z)-(f)_{y,r}|\,\d z\bigg]\cr
  &\leq &\tilde C_dC_d\big(M_R|\nabla f|(x)
  +M_R|\nabla f|(y)\big)r.
  \end{eqnarray}
Now \eqref{th.6.2-1} follows from the triangular inequality and
\eqref{III}, \eqref{IV}:
  \begin{eqnarray*}
  |f(x)-f(y)|&\leq&|f(x)-(f)_{x,r}|+|(f)_{x,r}-(f)_{y,r}|+|(f)_{y,r}-f(y)|\cr
  &\leq&2^{1+d}C_dM_R|\nabla f|(x)\,r+\tilde C_dC_d\big(M_R|\nabla f|(x)
  +M_R|\nabla f|(y)\big)r\cr
  &&+2^{1+d}C_dM_R|\nabla f|(y)\,r\cr
  &=& C_d(2^{1+d}+\tilde C_d)|x-y|\big(M_R|\nabla f|(x)
  +M_R|\nabla f|(y)\big).
  \end{eqnarray*}
We obtain \eqref{th.6.2-1}. \fin

\end{document}